\documentclass{article}

\usepackage{arxiv}

\usepackage[utf8]{inputenc} 
\usepackage[T1]{fontenc}    
\usepackage{hyperref}       
\usepackage{url}            
\usepackage{booktabs}       
\usepackage{amsfonts}       
\usepackage{nicefrac}       
\usepackage{microtype}      
\usepackage{lipsum}		
\usepackage{graphicx}
\usepackage{doi}
\usepackage{physics}

\newcommand{\pd}[2]{\frac{\partial #1}{\partial #2}}

\title{Solving engineering eigenvalue problems with neural networks using the Rayleigh quotient}

\date{} 					

\author{
    Conor Rowan \\
	Aerospace Engineering\\
	University of Colorado Boulder\\
	Boulder, CO 80309 \\
	\texttt{conor.rowan@colorado.edu} 
    \And
    Kurt Maute \\
	Aerospace Engineering\\
	University of Colorado Boulder\\
	Boulder, CO 80309 \\
	\texttt{kurt.maute@colorado.edu} 
    \And
    Alireza Doostan \\
	Aerospace Engineering\\
	University of Colorado Boulder\\
	Boulder, CO 80309 \\
	\texttt{alireza.doostan@colorado.edu} 
    \And
    John Evans \\
	Aerospace Engineering\\
	University of Colorado Boulder\\
	Boulder, CO 80309 \\
	\texttt{john.evans@colorado.edu} 
   }

\renewcommand{\headeright}{}
\renewcommand{\undertitle}{}
\renewcommand{\shorttitle}{}

\hypersetup{
pdftitle={A template for the arxiv style},
pdfsubject={q-bio.NC, q-bio.QM},
pdfauthor={David S.~Hippocampus, Elias D.~Striatum},
pdfkeywords={First keyword, Second keyword, More},
}

\begin{document}
\maketitle

\begin{abstract}
From characterizing the speed of a thermal system's response to computing natural modes of vibration, eigenvalue analysis is ubiquitous in engineering. In spite of this, eigenvalue problems have received relatively little treatment compared to standard forward and inverse problems in the physics-informed machine learning literature. In particular, neural network discretizations of solutions to eigenvalue problems have seen only a handful of studies. Owing to their nonlinearity, neural network discretizations prevent the conversion of the continuous eigenvalue differential equation into a standard discrete eigenvalue problem. In this setting, eigenvalue analysis requires more specialized techniques. Using a neural network discretization of the eigenfunction, we show that a variational form of the eigenvalue problem called the "Rayleigh quotient" in tandem with a Gram-Schmidt orthogonalization procedure is a particularly simple and robust approach to find the eigenvalues and their corresponding eigenfunctions. This method is shown to be useful for finding sets of harmonic functions on irregular domains, parametric and nonlinear eigenproblems, and high-dimensional eigenanalysis. We also discuss the utility of harmonic functions as a spectral basis for approximating solutions to partial differential equations. Through various examples from engineering mechanics, the combination of the Rayleigh quotient objective, Gram-Schmidt procedure, and the neural network discretization of the eigenfunction is shown to offer unique advantages for handling continuous eigenvalue problems.

\end{abstract}


\section{Introduction}

\paragraph{} The success of neural networks for diverse tasks such as regression, classification, generative modeling, and dimensionality reduction is well known \cite{jiang_densely_2022,lek_application_1996, zhao_review_2024, lorente_image_2021,  regenwetter_deep_2022, herron_latent_2024, he_deep_2021}. Another area of machine learning (ML) research with particular relevance to the engineering community is the use of neural networks to analyze physical systems governed by ordinary and partial differential equations. Whereas fully data-driven approaches have been deployed on data generated from these systems \cite{azizzadenesheli_neural_2024, bhaduri_stress_2022, lu_deeponet_2021, li_fourier_2021}, "physics-informed neural networks" (PINNs) incorporate governing differential equations in the training process. In this setting, the neural network often acts as a representation of the solution to a particular initial or boundary value problem, and the parameters of the neural network are tuned according to some physics-based criteria. These methods were first explored over two decades ago \cite{lagaris_artificial_1998}, but have only recently gained popularity owing to increased access to computational resources and sophisticated open-source ML software libraries. By proposing the basic infrastructure of discretizing the solution with a neural network, minimizing the squared error of the governing equation, and enforcing initial/boundary conditions with penalties, early works laid the foundation for an ever-increasing interest in ML methods for scientific computing \cite{raissi_physics-informed_2019, sirignano_dgm_2018}. Going forward, we will call the squared error objective commonly used in PINNs the "strong form loss." We note that in addition to offering an approach to a forward solution of a partial differential equation (PDE), PINNs can also be used to solve inverse problems. Inverse problems make use of data in order to estimate missing parameters in the governing equation. Since their introduction, PINNs methods for forward and inverse problems have been explored in the context of many different physical models such as elasticity \cite{cai_deep_2020, kag_physics-informed_2024}, contact mechanics \cite{sahin_solving_2024}, heat transfer \cite{madir_physics_2024}, and fluid mechanics \cite{jin_nsfnets_2021}.   

\paragraph{} Other researchers have taken the basic idea of PINNs---representing the solution to an ordinary or partial differential equation with a neural network---and extended it to different formulations of the governing equation. For example, it is well known that many PDEs of interest correspond to the minimum of a particular "variational energy" functional. Thus, instead of minimizing the strong form loss, the PDE solution can be approximated by minimizing the energy. This was first explored with the Deep Ritz method \cite{e_deep_2017}, and later extended to problems in linear elasticity \cite{liu_deep_2023}, hyperelasticity \cite{abueidda_deep_2022}, and fracture mechanics \cite{manav_phase-field_2024, ghaffari_motlagh_deep_2023}. In addition to energy minimization, weak form solutions of PDEs have been studied in the PINNs literature. A weak solution is obtained by requiring that the strong form residual is orthogonal to a user-defined basis. This approach has been explored in works such as \cite{kharazmi_variational_2019, kharazmi_hp-vpinns_2021, shang_deep_2022, khodayi-mehr_varnet_2019}. A benefit of using the variational energy or weak form system as the objective is that the order of spatial differentiation on the neural network representing the solution is lower. This simultaneously speeds up computations using automatic differentiation and expands the class of admissible solutions by reducing continuity requirements.

\paragraph{} Eigenvalue problems frequently show up in many areas of science and engineering. The eigenvectors of the covariance matrix are used in Principal Component Analysis to find directions of maximum variance \cite{makiewicz_principal_1992}. In the setting of random processes, the stationary distribution of a Markov Chain is governed by an eigenvalue problem \cite{seneta_computing_1980}, and the eigenfunctions of the covariance kernel are used to discretize continuous random processes using the Karhunen-Loeve expansion \cite{fontanella_karhunenloeve_2012}. In control theory, eigenvalues of the state-space system are used to characterize stability properties \cite{brogan_modern_1991}. In acoustics, the Helmholtz equation is an eigenvalue problem for the natural modes of vibration of the solution to the wave equation \cite{tai_helmholtzequation_1974}. In heat transfer, eigenvalues give insight into the decay rate of eigenfunction temperature distributions \cite{gerstenmaier_heat_2002}. In solid mechanics, eigenvalue problems are solved in order to determine the frequencies of vibration and corresponding mode shapes of a structure \cite{bathe_solution_1973} and are also used in buckling analysis \cite{bauchau_structural_2009}. 

\paragraph{} At first glance, eigenvalue differential equations appear to be a straightforward extension of PINNs methodologies. However, because both the eigenvalue and eigenfunction are unknown, and because solutions are non-unique, these problems present difficulties for out-of-the-box PINNs solvers. In fact, neural network discretizations of partial differential equations have seen limited treatment in the context of eigenvalue problems from engineering mechanics, despite being discussed in the original Deep Ritz Method paper \cite{e_deep_2017}. Physics-informed neural networks have been used to solve eigenvalue problems in quantum mechanics \cite{holliday_solving_2023, jin_physics-informed_2022} and for eigenfunctions of the Sturm-Liouville differential equation \cite{kovacs_conditional_2022, vigouroux_deep_2024}. Another work learns eigenfunctions for general self-adjoint operators \cite{ben-shaul_deep_2023}. In \cite{yang_neural_2023}, eigenanalysis for continuous eigenvalue differential equations is conducted with the power method using a neural network discretization of the eigenfunction. Though the power method is primarily for finding the first eigenvalue, it can be used with a "shift" operation to determine successive eigenvalues. However, this relies on prior knowledge of the distribution of eigenvalues in order to find the whole spectrum. This is an overly restrictive assumption in many problems of practical interest. To the best of the authors' knowledge, the only neural network discretization used for eigenvalue analysis specifically in the context of engineering mechanics appears in \cite{yoo_physics-informed_2024}, where the residual of the "strong form" of the eigenvalue problem with unknown eigenvalue is minimized. The authors report convergence issues in finding higher frequency eigensolutions with this formulation. This is due to the fact that the network tends to converge to the first eigenvalue and eigenfunction unless specifically prevented from doing so and that avoiding the lowest frequency eigenfunction is difficult to achieve by modifying the loss with a penalty. When using penalty methods---as opposed to strict constraint enforcement---the optimizer converges to local minima as opposed to finding successive eigenvalues/eigenfunctions. Noting the relative scarcity of previous work on eigenvalue analysis in engineering mechanics, and the issues posed by solving the eigenvalue problem with the power method, strong form error minimization, and penalty methods, our primary contributions in this work are as follows:

\begin{enumerate}
    \item We show how eigenvalue analysis with neural networks can be cast as a series of optimization problems using the Rayleigh quotient objective along with Gram-Schmidt orthogonalization;
    \item We demonstrate the efficacy of this approach by finding eigenvalues and eigenfunctions on a number of domains and show that the eigenfunctions are useful as a spectral basis;
    \item Taking advantage of the ability of neural networks to handle parameters, we use our method to perform uncertainty quantification and find a spectral basis on a parameterized geometry;
    \item We show that a neural network discretization outcompetes a standard spectral discretization on a nonlinear eigenvalue problem;
    \item We demonstrate the efficacy of the Rayleigh quotient method with neural network discretizations on a high-dimensional eigenvalue problem.
\end{enumerate}

The rest of this paper is organized as follows. In Section 2, we introduce necessary preliminaries such as the definition of the eigenproblem, a sketch of neural network discretizations of PDEs, and an overview of previous approaches to the eigenvalue problem taken in the literature. In Section 3, we introduce the Rayleigh quotient, Gram-Schmidt orthogonalization, and the iterative approach to eigenanalysis. In Section 4, we use this method to learn a harmonic basis on a variety of domains and show that the spectral basis can accurately approximate the solution to PDEs. In Section 5, we demonstrate that neural networks can easily be modified to solve parametric eigenvalue problems by performing uncertainty quantification on the lowest natural frequency of a multi-material linearly elastic structure and finding a spectral basis on a parameterized geometry. In Section 6, we test the Rayleigh quotient method on a nonlinear eigenvalue problem, and show that it is competitive with a more traditional approach. In Section 7, we study the Rayleigh quotient method for the $d$-dimensional Laplace equation and use Monte Carlo integration to approximate eigenvalues in high dimensions. Section 8 offers concluding thoughts and directions for future research.


\section{Preliminaries}

\subsection{Abstract Eigenvalue problem}

\paragraph{} The most general form of a continuous eigenvalue differential equation can be written as

\begin{equation}\label{generaleig}
    \begin{aligned}
        \mathcal{N}(\mathbf{u}_i)(\mathbf{x})  + \lambda_i h(\mathbf{u}_i)(\mathbf{x})  = \mathbf{0}, \quad \mathbf{x} \in \Omega, \\
        \mathcal{B}( \mathbf{u}_i)(\mathbf{x}) = \mathbf{0}, \quad \mathbf{x} \in \partial \Omega,
    \end{aligned}
\end{equation}

\noindent where $\mathcal{N}(\cdot)$ is a linear or nonlinear differential operator, $h(\cdot)$ is a function, $\mathbf{x}$ is the spatial coordinate, $\Omega$ is the domain, $\mathcal{B}(\cdot)$ is the boundary operator, $\lambda_i$ is the eigenvalue, and $\mathbf{u}_i(\mathbf{x})$ is the corresponding eigenfunction. Note that the eigenvalues are indexed by $i=1,2,\dots$ where $\lambda_1 < \lambda_2 < \dots$ by definition. We work with a nonlinear eigenvalue problem of the form given in Eq. \eqref{generaleig} in Section 6, however our emphasis is on so-called linear eigenvalue problems, which can be written as follows: 

\begin{equation*}\label{lineareig}
    \begin{aligned}
        \mathcal{L}(\mathbf{u}_i)(\mathbf{x}) + \lambda_i  \mathbf{u}_i(\mathbf{x}) = \mathbf{0}, \quad \mathbf{x} \in \Omega, \\
        \mathcal{B}( \mathbf{u}_i)(\mathbf{x}) = \mathbf{0}, \quad \mathbf{x} \in \partial \Omega,
    \end{aligned}
\end{equation*}

\noindent where $\mathcal{L}(\cdot)$ is a linear differential operator. Problems of this sort are used to find natural frequencies and modes of a linearly elastic structure, a case which we treat in Section 5. Note that the set of eigenvalues and eigenfunctions depends on the differential operator, the domain geometry, and the boundary conditions. Within the class of linear eigenvalue problems, we are particularly interested in the following Laplace eigenvalue problem

\begin{equation}\label{harmonic}
    \begin{aligned}
        \nabla^2 \mathbf{u}_i(\mathbf{x}) + \lambda_i  \mathbf{u}_i(\mathbf{x}) = \mathbf{0}, \quad \mathbf{x} \in \Omega, \\
        \mathcal{B}( \mathbf{u}_i)(\mathbf{x})  = \mathbf{0}, \quad \mathbf{x} \in \partial \Omega.
    \end{aligned}
\end{equation}

We will call functions which satisfy Eq. \eqref{harmonic} "harmonic" eigenfunctions. Noting that the Laplace operator is self-adjoint under the $L^2(\Omega)$ inner product, it follows that harmonic eigenfunctions corresponding to distinct eigenvalues are orthogonal with respect to that inner product.

\subsection{Neural network discretization}

\paragraph{} We use multilayer perceptron (MLP) neural networks to discretize the solution to eigenvalue problems of the sort discussed above. The neural network consists of the repeated application of the following transformation:

\begin{equation*}
    \mathbf{y}_k = \sigma\Big(  \mathbf{W}_k\mathbf{x}_k + \mathbf{B}_k  \Big),
\end{equation*}

\noindent where $\mathbf{x}_k$ is the input to the $k$-th layer, $\mathbf{W}_k$ are the weights, and $\mathbf{B}_k$ are the biases. The output $\mathbf{y}_k$ becomes the input to the next layer $\mathbf{y}_k \rightarrow \mathbf{x}_{k+1}$. The function $\sigma(\cdot)$ is an activation function that is applied element-wise. The trainable parameters of the neural network are grouped together as $\boldsymbol{\theta} = [ \mathbf{W}_1 , \mathbf{B}_1, \mathbf{W}_2,\mathbf{B}_2,\dots, \mathbf{W}_{\ell}]$ where $\ell$ is the number of layers in the network. The final layer is thought of as a dot product of spatial functions with coefficients, and thus does not have biases associated with it. Per \cite{sukumar_exact_2022, wang_exact_2023-1}, we build in the Dirichlet boundary conditions directly into the neural network discretization. The neural network discretization of the $i$-th eigenfunction can be written as

\begin{equation}\label{dirichlet}
\begin{aligned}
    \mathbf{\hat u}_i(\mathbf{x};\boldsymbol{\theta}) = \mathbf{D}(\mathbf{x}) \mathbf{\hat U}_i( \mathbf{x};\boldsymbol{\theta}) + \mathbf{G}(\mathbf{x}), \\
    \mathbf{D}(\mathbf{x}) = \mathbf{0}, \quad \mathbf{x} \in \partial \Omega_D, \\
    \mathbf{G}(\mathbf{x}) = \mathbf{g}(\mathbf{x}), \quad \mathbf{x} \in \partial \Omega_D.     
\end{aligned}
\end{equation}

The notation $\partial \Omega_D$ indicates the subset of the boundary with Dirichlet boundary conditions and $\mathbf{g}(\mathbf{x})$ gives the prescribed value along this part of the boundary. We call the discretized solution $\mathbf{\hat u}$ in order to disambiguate the exact solution of the eigenvalue problem from the approximate solution. The function $\mathbf{\hat U}(\mathbf{x};\boldsymbol \theta)$ is a neural network which need not satisfy the boundary conditions. As Eq. \eqref{dirichlet} shows, the Dirichlet boundaries are satisfied by construction. This contrasts with the Neumann boundary conditions, which must be enforced in the loss function. 

\paragraph{} For parametric problems, it is possible to approximate the dependence of the solution on the parameters by including them as additional inputs to the neural network. For a set of parameters $\mathbf{a}$---which might control the source term, boundary conditions, geometry, and/or coefficients of the PDE---the solution can be discretized in both physical and parameter space with

\begin{equation}\label{dirichletparam}
\begin{aligned}
    \mathbf{\hat u}_i(\mathbf{x},\mathbf{a};\boldsymbol{\theta}) = \mathbf{D}(\mathbf{x},\mathbf{a}) \mathbf{\hat U}_i( \mathbf{x},\mathbf{a};\boldsymbol{\theta}) + \mathbf{G}(\mathbf{x},\mathbf{a}), \\
    \mathbf{D}(\mathbf{x},\mathbf{a}) = \mathbf{0}, \quad \mathbf{x} \in \partial \Omega_D(\mathbf{a}), \\
    \mathbf{G}(\mathbf{x},\mathbf{a}) = \mathbf{g}(\mathbf{x},\mathbf{a}), \quad \mathbf{x} \in \partial \Omega_D(\mathbf{a}).     
\end{aligned}
\end{equation}

\subsection{Strong form loss of the eigenvalue problem}

\paragraph{} To simplify the exposition, we will work specifically with the eigenvalue problem given by Eq. \eqref{harmonic} to discuss various solution approaches. We note that if the solution is expanded as a linear combination of known spatial basis functions, this continuous eigenvalue problem can easily be transformed into a standard matrix eigenvalue problem. This is the typical approach when spectral or finite element bases are used in the solution of the eigenvalue problem \cite{vallet_spectral_2008}, but it is not possible in the case of neural network discretizations owing to the nonlinear dependence of the solution on the parameters $\boldsymbol \theta$. Inspired by the usual approach to PINNs-based solutions to PDEs, a first pass at finding $N$ eigenvalue-eigenfunction pairs is given by minimizing the strong form loss:

\begin{equation*}\label{firstpass}
    \underset{\boldsymbol{\theta},\boldsymbol{\lambda}}{\text{argmin }} \frac{1}{2} \sum_{i=1}^N \int_{\Omega} \Big \lVert 
 \nabla^2 \mathbf{\hat u}_i(\mathbf{x};\boldsymbol{\theta}) + \lambda_i \mathbf{\hat u}_i(\mathbf{x};\boldsymbol{\theta}) \Big \rVert^2 d\Omega,
\end{equation*}

\noindent where the neural network is structured such that there are $N$ outputs which each represent one of the eigenfunctions. The eigenvalues are introduced as additional learnable parameters. The problem with this approach is that nothing prevents the network from learning the first eigenpair $N$ times. In other words, there is no guarantee that each of the $N$ eigenvalues and eigenfunctions are distinct. Going forward, we will restrict our attention to self-adjoint operators so that eigenfunctions are orthogonal. In this setting, a penalty can be introduced which enforces orthogonality of the eigenfunctions to the residual of the eigenvalue equation. This is the approach taken in \cite{ben-shaul_deep_2023} to prevent redundant eigenpairs. In this case, the objective becomes

\begin{equation}\label{secondpass}
    \underset{\boldsymbol{\theta},\boldsymbol{\lambda}}{\text{argmin }} \frac{1}{2} \sum_{i=1}^N \int_{\Omega} \Big \lVert 
 \nabla^2 \mathbf{\hat u}_i(\mathbf{x};\boldsymbol{\theta}) + \lambda_i \mathbf{\hat u}_i(\mathbf{x};\boldsymbol{\theta}) \Big \rVert^2 d\Omega + \frac{\beta}{2} \sum_{i=1}^N\sum_{j=1}^i \qty(\int_{\Omega} \mathbf{\hat u}_i(\mathbf{x};\boldsymbol{\theta}) \cdot \mathbf{\hat u}_j(\mathbf{x};\boldsymbol{\theta}) d\Omega - \delta_{ij})^2 ,
\end{equation}

\noindent where $\beta$ is a penalty parameter and $\delta_{ij}$ is the Kronecker delta symbol. 

\paragraph{} Say that we take $\mathbf{\hat u}(\mathbf{x})$ in Eq. \eqref{secondpass} to be a scalar on the 1D domain $x\in[0,1]$ with homogeneous Dirichlet boundaries. We can gauge the feasibility of this objective by eliminating the first term which penalizes error in the eigenvalue equation. In other words, we want to learn an orthonormal set of functions which need not be harmonic eigenfunctions. The objective here is simply

\begin{equation}\label{L1}
\mathcal{L}_1 = \frac{1}{2}\sum_{i=1}^N\sum_{j=1}^i \qty(\int_0^1 \hat u_i(x;\boldsymbol{\theta}) \hat u_j(x;\boldsymbol{\theta}) - \delta_{ij})^2.
\end{equation}

\begin{figure}[hbt!]
\centering
\includegraphics[width=0.65\textwidth]{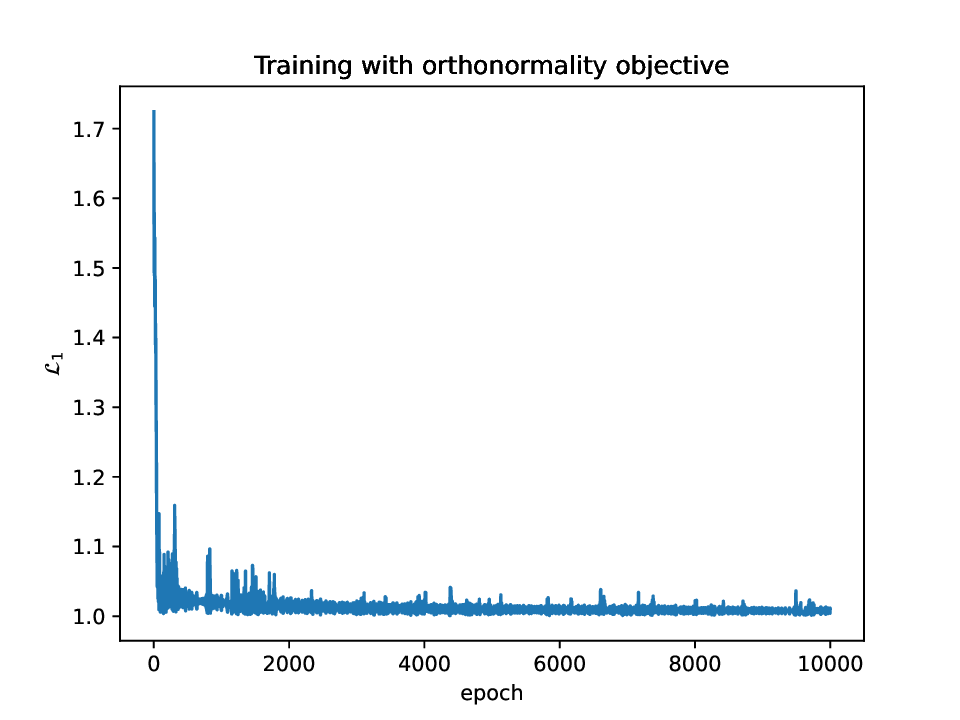}
\caption{Even when the eigenvalue differential equation is ignored in Eq. \eqref{secondpass} and $N=3$ in one spatial dimension, the optimization problem is intractable. The network should obtain $\mathcal{L}_1=0$ when the functions are orthogonal, but remains stuck in a local minimum.}
\label{pp1}
\end{figure}

Taking $N=3$, we minimize the orthogonality penalty in terms of the neural network parameters which construct the three functions. See Figure \ref{pp1} for a representative training trajectory. Even for this simple problem, the optimizer spends thousands of epochs stuck in a local minimum. This stagnation of the training occurs even though the network is sufficiently expressive to accurately represent three orthogonal functions if explicitly trained to approximate them. Due to tradeoffs between the normalization condition and orthogonality, this penalty is a complex loss function even in the absence of the Laplacian eigenproblem residual. That being said, we observe similar stagnation of the loss when using the objective given by Eq. \eqref{secondpass}. Thus, in this example, incorporating the residual of the differential equation does not provide any regularization to the optimization problem. In \cite{yoo_physics-informed_2024}, instead of using the orthogonality of the eigenfunctions, redundant eigenpairs are avoided by enforcing that $\lambda_{i+1}>\lambda_i$ with a penalty. Numerical experimentation indicates that this also leads to complex optimization problems that tend to get stuck in local minima, even for 1D scalar eigenvalue problems. These convergence issues are noted by the authors in their paper.

\paragraph{} As noted in the original Deep Ritz paper \cite{e_deep_2017}, there is an alternative variational formulation of eigenvalue problems called the "Rayleigh quotient." The Rayleigh quotient objective reduces the order of differentiation on the solution and removes dependence of the objective on the unknown eigenvalue. Removal of the eigenvalue as an additional unknown in the optimization problem avoids problems with non-uniqueness as noted in \cite{yoo_physics-informed_2024}. Furthermore, the Rayleigh quotient is naturally phrased as an optimization problem and speeds up automatic differentiation calculations by lowering the order of differentiation on the eigenfunctions. These considerations suggest that the Rayleigh quotient is much better suited for neural network discretizations of eigenvalue problems. We note that the Rayleigh quotient is used as the objective in \cite{ben-shaul_deep_2023} (instead of the strong form loss) but with the orthogonality penalty given by Eq. \eqref{L1} to enforce the distinctness of the eigenfunctions. As this brief discussion shows, the orthogonality penalty introduces convergence issues for the optimization problem. Our goal is to obtain the ordered sequence of eigenpairs $\{ \lambda_i, \mathbf{\hat u}_i\}_{i=1}^N$ with $N$ as large as possible. As we discuss in the sequel, the Rayleigh quotient objective also naturally enforces that successive eigenpairs are ordered. Thus, given its numerical benefits (lower order of differentiation and fewer unknown optimization variables) and its ability to find ordered eigenpairs, we forego the strong form loss of the eigenvalue problem. Instead, we search for a way to use the Rayleigh quotient objective without relying on the penalty on the orthogonality of the eigenfunctions.


\section{Rayleigh quotient}

\paragraph{} As shown in \cite{saad_numerical_2011}, and as discussed in \cite{e_deep_2017, ben-shaul_deep_2023}, the Rayleigh quotient is an alternative formulation of the eigenvalue problem. For the Laplace eigenvalue problem given by Eq. \eqref{harmonic}, the first eigenvalue is found with

\begin{equation*}
    \lambda_1 = \underset{u_1(\mathbf{x})}{\text{min }} \frac{ \int_{\Omega} \nabla u_1 \cdot \nabla u_1 d\Omega }{\int_{\Omega} u_1^2 d\Omega},
\end{equation*}

\noindent where the minimizer $u_1(\mathbf{x})$ is the corresponding eigenfunction. Analogous Rayleigh quotients can be formed for structural mechanics problems with a given constitutive tensor. We note that the Rayleigh quotient enforces zero Neumann boundary conditions for boundaries which do not have a prescribed displacement per Eq. \eqref{dirichlet}. Successive eigenvalues and eigenfunctions can be found with

\begin{equation}\label{lam1}
    \lambda_i = \underset{u^{\perp}_i(\mathbf{x})}{\text{min }} \frac{ \int_{\Omega} \nabla u^{\perp}_i \cdot \nabla u^{\perp}_i d\Omega }{\int_{\Omega} (u_i^{\perp})^2 d\Omega}, \quad i=1,2,\dots
\end{equation}

The superscript "$\perp$" is used to indicate that the $i$-th eigenfunction is orthogonal to all the previous ones, i.e. $\int_{\Omega} u^{\perp}_i u_j d\Omega = 0 $ for all $j<i$. The $i$-th Rayleigh quotient, given by Eq. \eqref{lam1}, is minimized in a subspace which is orthogonal to all previous eigenfunctions. Whereas a penalty was used to satisfy this condition in previous works, we enforce orthogonality automatically by using the Gram-Schmidt process to orthogonalize the search space. Consider a discretized function given by $\tilde u_i(\mathbf{x};\boldsymbol{\theta})$. The following is orthogonal to the previous eigenfunctions by construction:

\begin{equation*}\label{gs}
    \hat u_i^{\perp}(\mathbf{x};\boldsymbol{\theta}) = \tilde u_i(\mathbf{x};\boldsymbol{\theta}) - \sum_{j=1}^{i-1} \qty(\frac{\int_{\Omega} \tilde u_i \hat u_j d\Omega}{\int_{\Omega} \hat u_j^2 d\Omega}) \hat u_j(\mathbf{x}).
\end{equation*}

The functions $\hat u_j$ are the previous converged approximations of  eigenfunctions. This means that

\begin{equation*}\label{lam}
    \hat u_i(\mathbf{x}) := \underset{\hat u^{\perp}_i(\mathbf{x})}{\text{argmin }} \frac{ \int_{\Omega} \nabla \hat u^{\perp}_i \cdot \nabla \hat u^{\perp}_i d\Omega }{\int_{\Omega} (\hat u_i^{\perp})^2 d\Omega}, \quad i=1,2,\dots
\end{equation*}

In other words, we identify the $i$-th eigenfunction as the minimizer of the Rayleigh quotient in the space of functions orthogonal to all previous eigenfunctions. The Gram-Schmidt process builds in the orthogonality restriction without a penalty. To obtain eigenfunctions with unit magnitude, we can append a penalty to the Rayleigh quotient to enforce this:

\begin{equation*}\label{lampen}
     \underset{\hat u^{\perp}_i(\mathbf{x})}{\text{argmin }} \frac{ \int_{\Omega} \nabla \hat u^{\perp}_i \cdot \nabla \hat u^{\perp}_i d\Omega }{\int_{\Omega} (\hat u_i^{\perp})^2 d\Omega} + \beta \qty(\int_{\Omega} (\hat u^{\perp}_i)^2 d\Omega - 1)^2, \quad i=1,2,\dots
\end{equation*}

The second term penalizes deviations of the eigenfunctions from unit magnitude. Penalties of this sort are used widely in the machine learning literature and thus we do not expect convergence issues like the orthogonality penalty of Eq. \eqref{L1}. Going forward, we depart from the Laplace eigenvalue problem and denote a generic Rayleigh quotient with a potentially vector-valued solution as $\mathcal{R}\Big(\mathbf{u}(\mathbf{x})\Big)$.

\paragraph{} As Eq. \eqref{lam1} shows, the eigenpairs are found by solving successive optimization problems. We now discuss the implementation of these optimization problems and the Gram-Schmidt procedure. We use a single neural network to discretize each of the $N$ eigenfunctions, whose parameters are overwritten in each successive optimization problem. The neural network output for the $i$-th "auxiliary function" is denoted $\mathbf{\tilde u}_i$. We call it an auxiliary function to be clear that it only becomes an eigenfunction once orthogonalized. The eigenvalues (the value of the Rayleigh quotient at its minimum) and the eigenfunctions (the minimizers of the Rayleigh quotient) will be stored on an integration grid once each optimization problem converges. The first $N$ eigenfunctions are found by solving the following optimization problem for $i=1,2,\dots,N$. The $i$-th optimization problem is

\begin{equation}\label{algorithm}
    \begin{aligned}
        \mathcal{L}^R_i = \mathcal{R} \Big( \mathbf{\hat u}_i^{\perp}(\mathbf{x};\boldsymbol{\theta}) \Big) + \beta \qty( \int_{\Omega} \lVert \mathbf{\hat u}^{\perp}(\mathbf{x};\boldsymbol{\theta}) \rVert^2 d \Omega - 1)^2,\\
        \mathbf{\hat u}_i^{\perp}(\mathbf{x};\boldsymbol{\theta}) = \mathbf{\tilde u}_i(\mathbf{x};\boldsymbol{\theta}) - \sum_{j=1}^{i-1} \qty(\frac{\int_{\Omega}  \mathbf{\tilde u}_i \cdot \mathbf{\hat u}_j d\Omega}{\int_{\Omega} \mathbf{\hat u}_j \cdot \mathbf{\hat u}_j d\Omega}) \mathbf{\hat u}_j(\mathbf{x}),\\
        \boldsymbol{\theta}_i = \underset{\boldsymbol{\theta}}{\text{argmin }} \mathcal{L}_i^R(\boldsymbol{\theta}).
    \end{aligned}
\end{equation}



If it is not of interest that the neural network directly outputs normalized eigenfunctions, we can remove the penalty on their magnitude by setting $\beta=0$. We will experiment with both approaches in the forthcoming numerical examples. Note that there is no dependence of the previous converged eigenfunctions on the current neural network parameters. Eq. \eqref{algorithm} can be implemented in PyTorch, where the Rayleigh quotient is formed with the orthogonalized candidate eigenfunction and spatial gradients are computed using automatic differentiation. The neural network parameters are determined such that the Rayleigh quotient is minimized using ADAM optimization. For all following examples, we use ADAM optimization with a fixed learning rate of $5 \times 10^{-3}$ unless otherwise noted to minimize the Rayleigh quotient. We note from experience that for Rayleigh quotients corresponding to higher frequency eigenvalues, ADAM requires more epochs to converge. Thus, ADAM optimization continues until the current Rayleigh quotient falls within a user-specified distance of the previous converged value, and then a fixed number of epochs are run. For all examples, a fixed number of epochs are run once $\lambda_i \leq \lambda_{i-1} + \mathcal{T}$ where $\mathcal{T}$ is a convergence threshold parameter chosen for each numerical example. We tend to increase this fixed number of epochs with the index of the eigenproblem, which can be seen in the convergence plots. Once a minimum is obtained, the current eigenfunction is returned and stored for use in orthogonalizing the next candidate eigenfunction.

\paragraph{} We note that the Gram-Schmidt procedure may also be used with the strong form loss in order to avoid penalties on the orthogonality of eigenfunctions or on the size of the eigenvalues. But, because the Rayleigh quotient at the minimum is the value of the eigenvalue, finding the global minima of the successive optimization problems in Eq. \eqref{algorithm} corresponds to finding the eigenvalues in ascending order. We have no guarantee that the optimizer will find such global minima, but we will observe whether the proper ordering of the eigenvalues is obtained in the numerical examples. In contrast, the strong form loss with the Gram-Schmidt procedure has no such property encouraging the ordering of the eigenvalues. Any eigenpair that satisfies the differential equation and is distinct from the previous eigenpairs is an equally valid solution. We believe this encouraging of the ordering of eigenvalues, as well as the fact that we avoid introducing the unknown eigenvalue into the optimization problem, to be arguments in favor of the Rayleigh quotient objective. We now investigate the efficacy of this method for various eigenvalue problems relevant to engineering mechanics.


\section{Learning spectral bases}

\paragraph{} Functions that satisfy Eq. \eqref{harmonic} are often used as a spectral basis to approximate smooth functions. Harmonic eigenfunctions of this sort are especially useful for discretizing smooth solutions to PDEs because they are orthogonal with respect to the $L^2$ inner product, encode the geometry and boundary conditions of the problem, and, in comparison to finite difference and finite element methods, can significantly reduce the number of degrees of freedom in the approximation while retaining accuracy. As such, obtaining spectral bases provides a strong motivation for the importance of obtaining a large set of eigenfunctions. Our first test case is finding 1D harmonic functions, for which we have exact solutions.

\subsection{1D Fourier basis}

\paragraph{} In one spatial dimension, Eq. \eqref{harmonic} with homogeneous Dirichlet boundaries reduces to

\begin{equation*}
    \frac{\partial^2 u_i}{\partial x^2} + \lambda_i u_i = 0, \quad u(0)=u(1)=0.
\end{equation*}

\begin{figure}[hbt!]
\centering
\includegraphics[width=0.65\textwidth]{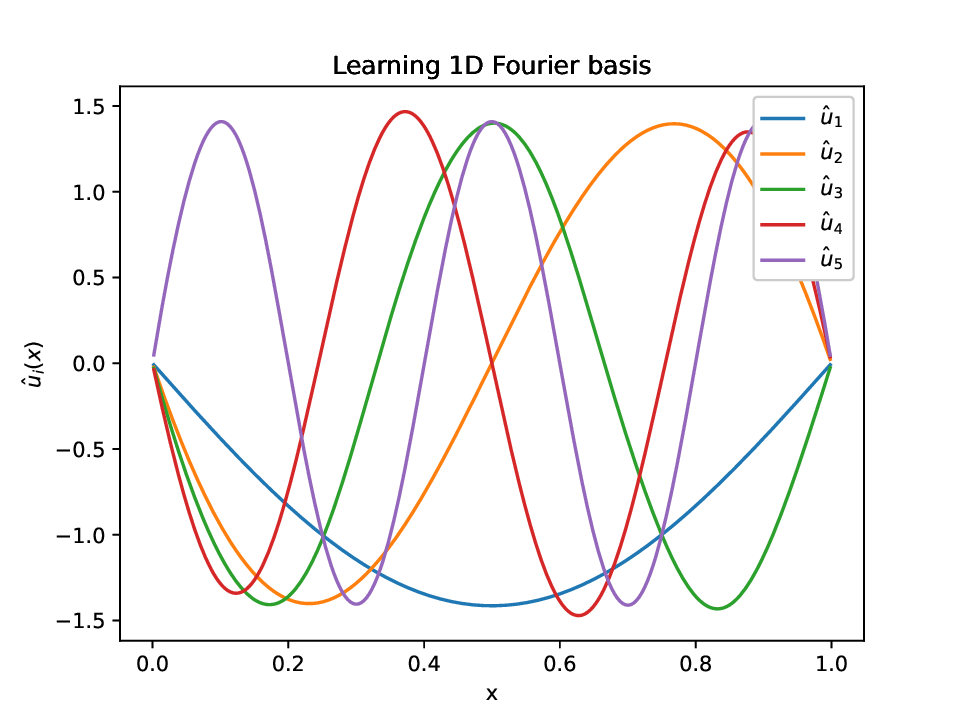}
\caption{Learning the orthonormal 1D Fourier basis with homogeneous Dirichlet boundaries. Though the penalty on the Rayleigh quotient enforces that the eigenfunctions have unit magnitude, they are only defined up to a minus sign.}
\label{pp2}
\end{figure}

The orthonormal set of eigenfunctions can be written down by inspection as $u_i(x)=\sqrt{2} \sin( i \pi x)$ with corresponding eigenvalues $\lambda_i = i^2 \pi^2 $ for $i=1,2,\dots$. We can verify that the proposed method recovers these functions. We will build in the boundary conditions to the neural network discretization of the basis functions, and use a single network for each optimization problem in Eq. \eqref{algorithm}, whose parameters are overwritten after storing the converged solution. We use a single hidden layer network with a width of 20 to learn the basis and hyperbolic tangent activation functions. Because the eigenvalues increase so rapidly in 1D, we take the convergence threshold parameter to be $\mathcal{T}=300$. The integration grid consists of 250 equally spaced points. In Figure \ref{pp2}, the first five basis functions of the orthonormal Fourier basis learned from Eq. \eqref{algorithm} are shown. In order to assess the accuracy of the learned basis, we define the following error metrics:

\begin{equation} \label{toterror}
    \mathcal{E}_u = \frac{1}{N}\sum_{i=1}^N \int_0^1 \Big( \hat u_i(x) - u_i(x) \Big)^2 dx, \quad \mathcal{E}_{\lambda} = \frac{1}{N}\sum_{i=1}^N\frac{|\hat \lambda_i - \lambda_i |}{\lambda_i},
\end{equation}

\noindent where $N$ is the number of terms in the learned basis. We understand the error in the eigenfunctions to be a relative error measure because they are orthonormal by construction. Our method can go up to $N=14$ before the optimizer converges to spurious solutions. We speculate that the causes of this are the spectral bias of neural networks, whereby high-frequency solutions are more difficult to represent \cite{rahaman_spectral_2019}, and the accumulation of numerical errors from the Gram-Schmidt process. The two error measures for $N=14$ are $\mathcal{E}_u=3\times 10^{-4}$ and $\mathcal{E}_{\lambda}=4 \times 10^{-3}$. The one-dimensional Fourier basis is a good first test case because simple analytical solutions exist to compare the learned harmonic functions against. We note that eigenfunctions become increasingly difficult to recover as the frequency increases.

\subsection{Semicircular domain}

\paragraph{} We now compute a harmonic set of functions on an irregular domain and show that they are an effective spectral basis for approximating solutions to PDEs. Like neural network discretizations, spectral bases are of interest because they obviate the need for meshing the computational domain. Unlike neural networks, they do not introduce additional nonlinearity to the problem, meaning that they can dramatically reduce the computational cost of the PDE solution. We take the computational domain $\Omega$ to be the unit semicircle. The domain has homogeneous Dirichlet boundary conditions, and we solve Eq. \eqref{harmonic} with the Rayleigh quotient method given in Eq. \eqref{algorithm}. We use a two hidden-layer neural network with a width of 10 and sigmoid activation functions in the following. The convergence threshold parameter is $\mathcal{T} = 25$ given that the eigenvalues increase less rapidly in higher dimensions. We note that it is not necessary to use very large and/or deep neural networks to represent these functions. The boundary conditions are enforced per Eq. \eqref{dirichlet} by multiplying by the neural network output by $D(x_1,x_2) = x_2 ( 1 - x_1^2 - x_2^2 )$, which is a distance level-set type function. We note from experience that better convergence is observed when there is no penalty used in Eq. \eqref{algorithm}, and instead the eigenfunctions are normalized once the Rayleigh quotient converges. In the 1D problem, the penalty seemed to present fewer numerical issues than the 2D semicircular problem, for which the penalty parameter causes the Rayleigh quotient to stagnate in local minima. 

\begin{figure}[hbt!]
\centering
\includegraphics[width=1.1\textwidth]{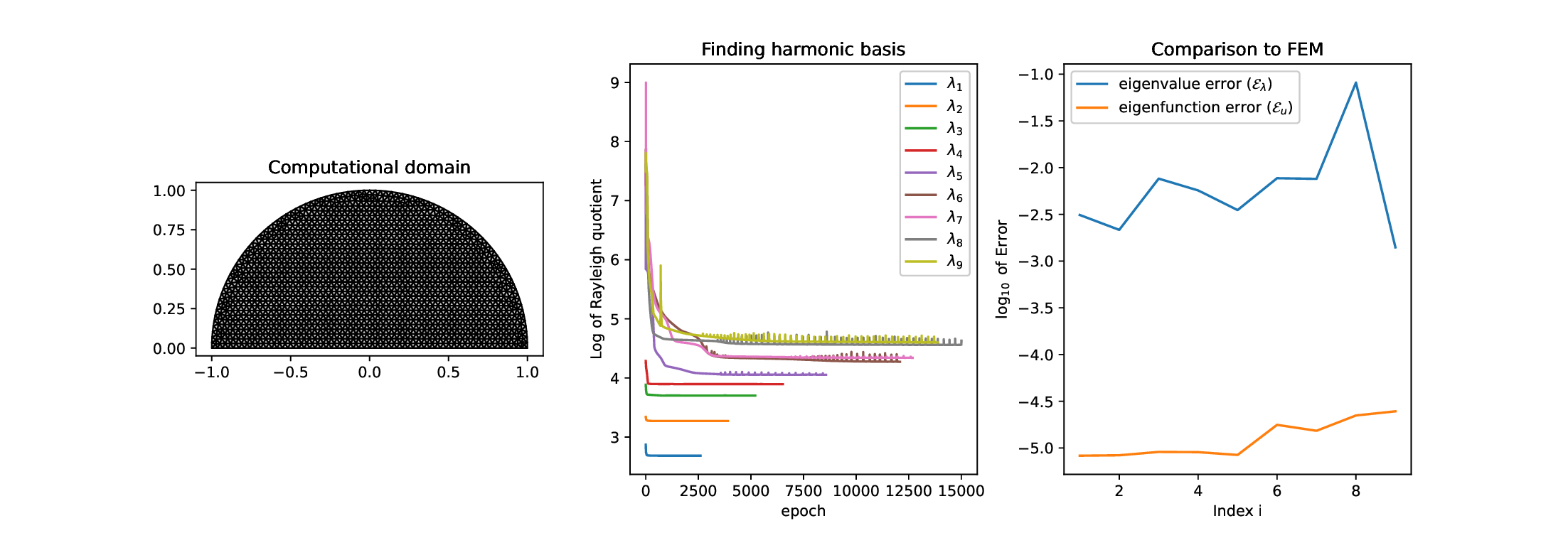}
\caption{The semicircular domain on which we find the harmonic basis and the mesh used in the finite element solution (left). The minimization of the Rayleigh quotient without penalty carries on for a fixed number of epochs once it falls below a threshold set by the previous converged value (center). The error with the finite element method results for both the eigenfunctions and eigenvalues (right). We note that the Rayleigh quotient objective provides ordered eigenpairs.}
\label{domain}
\end{figure}

\begin{figure}[hbt!]
\centering
\includegraphics[width=0.85\textwidth]{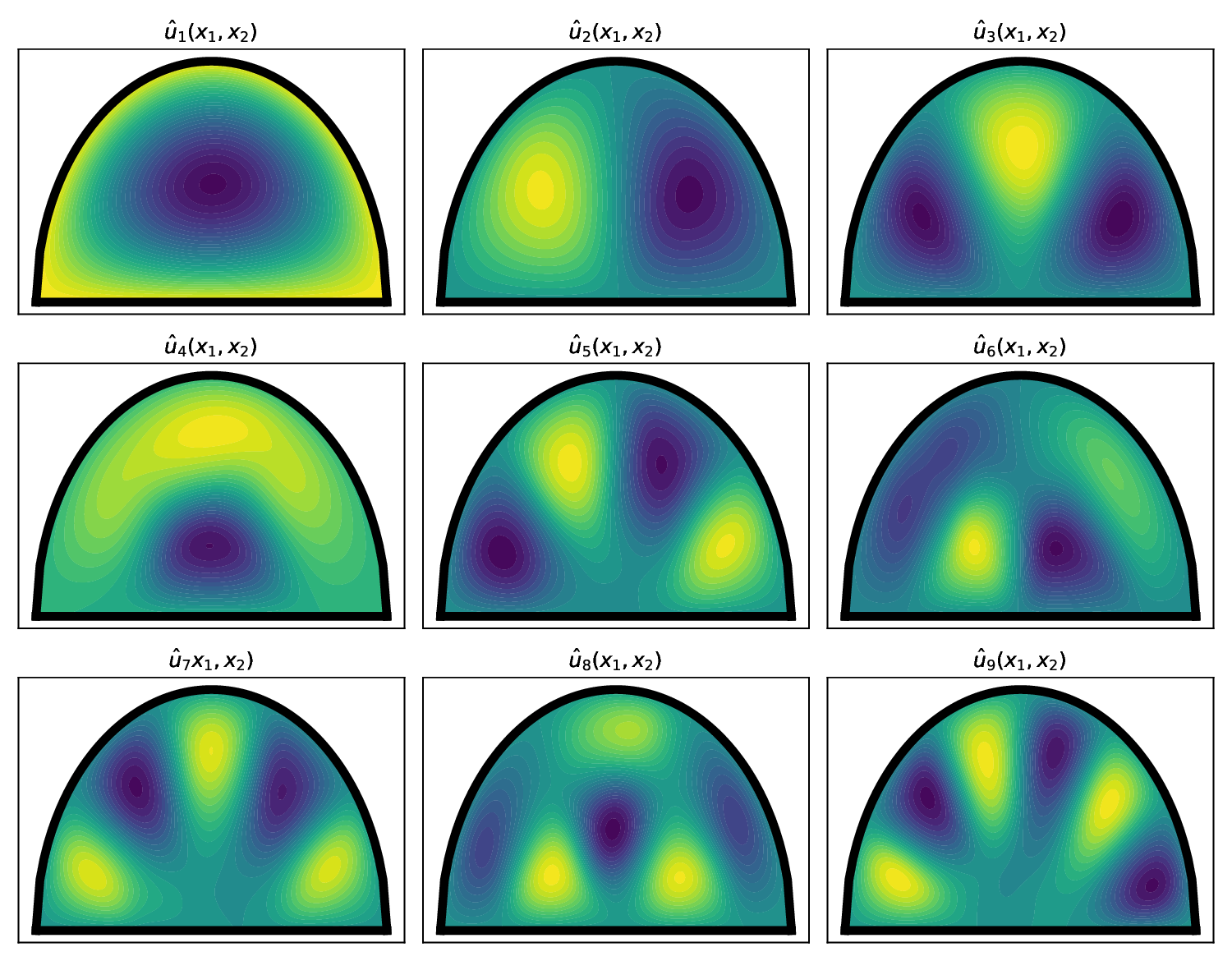}
\caption{The first 9 harmonic basis functions on the semicircular domain with homogeneous Dirichlet boundaries. Given that we do not enforce a penalty on the magnitude of the eigenfunctions, their magnitude is not uniquely defined by the Rayleigh quotient. As both the sign and the magnitude of the eigenfunctions are arbitrary, we do not include a colorbar in the visualization.}
\label{semicirclebasis}
\end{figure}

\paragraph{} We compare our neural network discretization of the eigenproblem to results obtained with the finite element method (FEM). The Python implementation of the FEM toolbox CALFEM is used to mesh the semicircular domain and solve the eigenvalue problem \cite{ottosson_implementation_2010}. With FEM, eigenanalysis is carried out by forming the mass and stiffness matrices with triangular elements and a linear interpolation of the solution field. We then use standard techniques to solve matrix eigenvalue problems for the eigenvalues and eigenfunctions stored on the nodes of the mesh. The two error metrics for eigenfunction and eigenvalue will be computed for each eigenproblem as 

\begin{equation*}
    \mathcal{E}_{u_i} =  \int_{\Omega} \Big( \hat u_i(\mathbf{x}) - u^{\text{FEM}}_i(\mathbf{x}) \Big)^2 d\Omega, \quad \mathcal{E}_{\lambda_i} = \frac{|\hat \lambda_i - \lambda^{\text{FEM}}_i |}{\lambda^{\text{FEM}}_i},
\end{equation*}

\noindent where $i=1,2,\dots,N$ is the index of the eigenproblem. We note that the eigenfunctions from both methods are normalized, so there is no need to include a denominator to obtain a relative error measure. In the finite element solution, we use a mesh with 3060 elements and verify that it is sufficiently refined by computing the percent change of the highest eigenvalue of interest with further refinement. In this example, we are interested in the first $N=9$ eigenpairs. When the mesh is refined to 4153 elements, $\lambda_9$ changes by $0.2\%$, which we deem to be acceptable convergence. See Figure \ref{domain} for the computational domain, the training convergence of the first 9 Rayleigh quotients, and the error of the computed eigenvalue and eigenfunction with the FEM results. The maximum error for the eigenvalues is $\text{max}(\mathcal{E}_{\lambda_i})=1.5\times10^{-2}$ and the maximum squared error for the eigenfunctions is $\text{max}(\mathcal{E}_{u_i})=7 \times 10^{-5}$. Figure \ref{semicirclebasis} shows the corresponding eigenfunctions. As expected, the eigenvalues---understood as the converged value of the Rayleigh quotient---are ordered. Higher index eigenfunctions correspond to higher frequency behavior, which is what we expect from a harmonic basis.

\paragraph{} We now use the spectral basis to discretize solutions to PDEs. We note that in a recent work \cite{hao_laplacian_2025}, the eigenfunctions of the Laplacian operator are used in the context of neural operators in order to discretize a missing nonlinearity in the governing partial differential equation. As they show, this approach demonstrates faster convergence than when the basis is learned from data with a neural network. However, in their work, the eigenfunctions are obtained via conventional numerical techniques such as the finite element method. In part, we build on this work by proposing an alternative neural network-based technique to obtain the eigenfunctions, which is convenient if one is already doing operator learning with machine learning libraries. In our case, instead of discretizing missing physics, we solve the governing equation for steady state heat conduction with isotropic material using the eigenfunction basis. The temperature $T(x_1,x_2)$ is governed by 

\begin{equation}\label{manufactured}
\begin{aligned}
    \nabla^2 T + s(\mathbf{x}) = 0, \quad \mathbf{x} \in \Omega, \\
    T(\mathbf{x}) = 0, \quad \mathbf{x} \in \partial \Omega.
\end{aligned}
\end{equation}

The source term $s(\mathbf{x})$ is determined by the method of manufactured solutions. We assume a parameterized form of the manufactured solution for the temperature, compute the corresponding source term, use this source term with the spectral discretization to solve the PDE, then compare the approximate solution to the exact one. The temperature is parameterized by 

\begin{equation*}
    T(x_1,x_2) = D(\mathbf{x})\sum_{0 \leq(i+j)\leq 3} t_{ij} x_1^ i x_2^j,
\end{equation*}

\noindent where each component of the parameter $\mathbf{t}$ is independent and uniformly distributed with $t_{ij} \sim \mathcal{U}(0,3)$. The function $D(\mathbf{x})$ enforces the homogeneous Dirichlet boundary conditions. Assuming a solution governed by stochastic parameters is one way to assess the accuracy of our spectral basis on a certain class of problems. The approximate solution is discretized with the eigenfunctions as $\hat T(\mathbf{x}) = \sum_{i=1}^N a_i \hat u_i(\mathbf{x})$. The Galerkin weak form solution for the coefficients on the harmonic basis functions for $j$-th realization of the random temperature distribution is 

\begin{equation*}\label{galerkin}
    \begin{aligned}
        \mathbf{a}^j = \mathbf{K}^{-1} \mathbf{F}^j, \\
        K_{k\ell}= \int_{\Omega} \nabla \hat u_k \cdot \nabla \hat u_{\ell} d\Omega, \\
        F^j_k = \int_{\Omega} \Big( -\nabla^2 T_j \Big) \hat u_k d\Omega.
    \end{aligned}
\end{equation*}

The index $j$ denotes a random realization of the parameters controlling the manufactured solution and we take $S$ samples total. We compute the average normalized error $\mathcal{E}$ between the approximated and real solution as 

\begin{equation}\label{avgerror1}
    \mathbb{E}\Big( \mathcal{E}(\hat T - T) \Big) \approx \frac{1}{S}\sum_{j=1}^S \frac{\int_{\Omega} \Big( T_j - \sum_{i=1}^N a_i^j \hat u_i(\mathbf{x}) \Big)^2 d\Omega}{\int_{\Omega} \lVert T_j \rVert^2 d\Omega}.
\end{equation}

Using $N=15 $ harmonic basis functions learned as described above, we compute the approximate solution for $S=10000$ samples. The error is normalized by the squared magnitude of the exact solution. The average error we obtain is $\mathbb{E}\Big( \mathcal{E}(\hat T - T ) \Big)=7\times10^{-3}$, which is small given that Eq. \eqref{avgerror1} is a relative error metric. Thus, the spectral basis performs well for the polynomial-type forcing. We note that for higher frequency solutions, the number of terms in the spectral basis needed for high accuracy quickly grows. In addition to showcasing the method for obtaining the eigenpairs, our intention here is to suggest one possible use case of the eigenfunctions.

\section{Parametric eigenvalue problems}

\subsection{Uncertainty quantification}

\paragraph{} We modify the neural network discretization to solve a family of parameterized eigenvalue problems simultaneously. This is often useful in the context of uncertainty quantification, as we will demonstrate. One interesting parametric eigenvalue problem is finding the first natural frequency of a structure with spatially varying material properties. Using index notation, the parametric eigenvalue problem for structural dynamics assuming linear elastic constitutive behavior is given by

\begin{equation*}
    \pd{}{x_j}C_{ijk\ell}(\mathbf{x};\mathbf{a}) \pd{u_k}{x_{\ell}} + \lambda u_i = 0,
\end{equation*}

\noindent where the vector $\mathbf{a}$ collects the parameters that control the spatial distribution of material in the structure. Using Voigt notation, the Rayleigh quotient corresponding to this eigenproblem is 


\begin{equation}\label{param1}
    \lambda(\mathbf{a}) = \underset{ \mathbf{u}(\mathbf{x};\mathbf{a})}{\text{min }} \frac{\int_{\Omega} \tilde C_{ij}(\mathbf{x};\mathbf{a}) \epsilon_i \epsilon_j d \Omega }{ \int_{\Omega} u_i u_i d\Omega},
\end{equation}

 \noindent where the strain vector is given by $\boldsymbol \epsilon = [ \partial u_1 / \partial x_1,\partial u_2/\partial x_2, \partial u_1/\partial x_2 + \partial u_2 / \partial x_1]$. Assuming the structure is in a state of plane stress, the constitutive matrix is

\begin{equation}\label{const}
    \mathbf{\tilde C}(\mathbf{x};\mathbf{a}) = \frac{E(\mathbf{x};\mathbf{a})}{1-\nu^2} \begin{bmatrix} 1 & \nu & 0 \\ \nu & 1 & 0 \\ 0 & 0 & (1-\nu)/2
    \end{bmatrix},
\end{equation}

\noindent where $\nu$ is the Poisson ratio and $E(\mathbf{x};\mathbf{a})$ is the parameterized modulus. An example of a parameterization of this sort is a multimaterial structure with the location of the interface between the two materials controlled by a single random parameter:

\begin{equation}\label{material}
    E(x_1,x_2;a) = \frac{E_1-E_0}{2} \Big( \tanh\Big(q(x_1-a)\Big) \Big) + E_0.
\end{equation}

\begin{figure}[hbt!]
\centering
\includegraphics[width=0.7\textwidth]{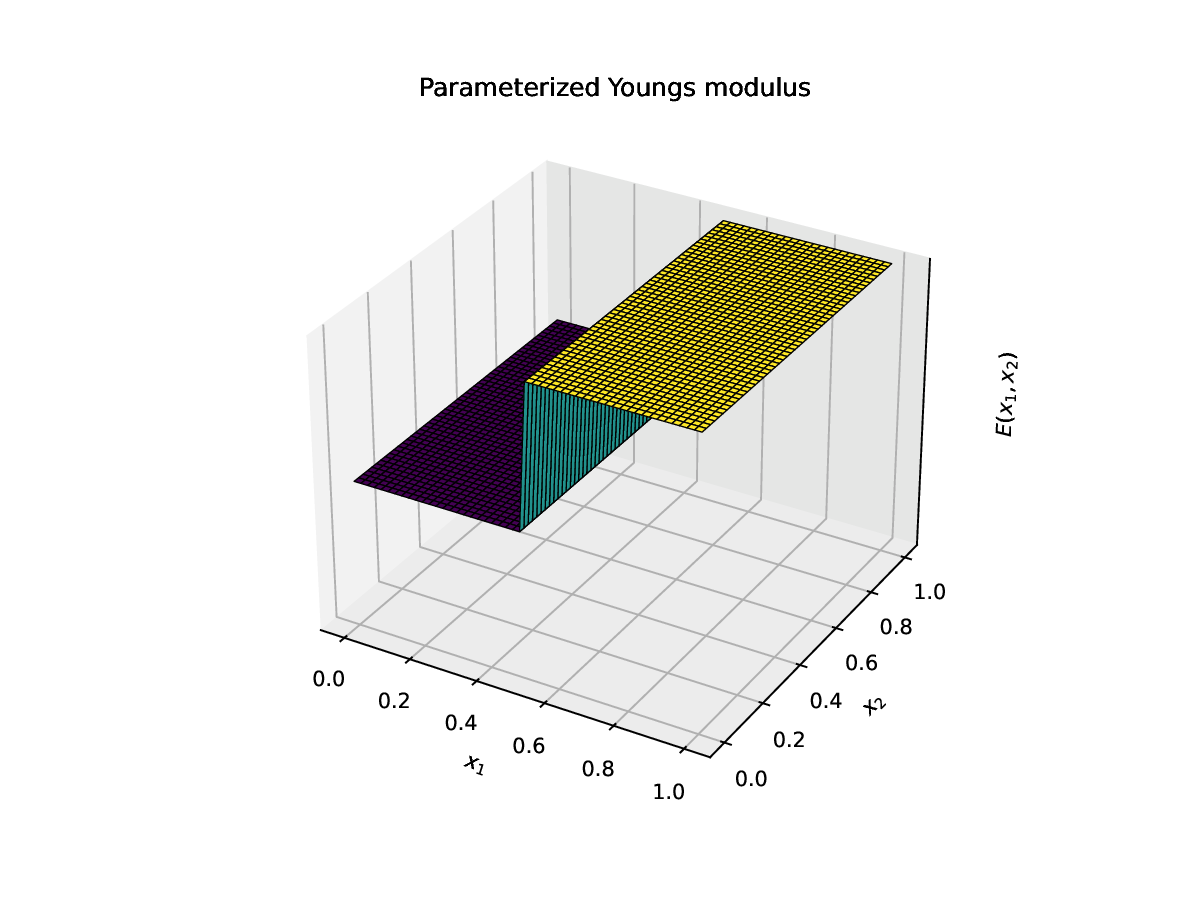}
\caption{Parameterized Young's modulus for the multimaterial structure. The interface between the two materials is at $x_1=a$. In the figure, the interface position is $x_1=a=1/2$. }
\label{diagram}
\end{figure}

When the parameter $q$ is chosen to be large, this expression accurately approximates a jump in the modulus from $E_0$ to $E_1$ with an interface at $x_1=a$. See Figure \ref{diagram} for a schematic of the problem set up. 

\paragraph{} As shown in Eq. \eqref{dirichletparam}, the parameters are introduced as additional inputs to the neural network. The neural network representation of the two displacement components is then given by $\mathbf{\hat u}(\mathbf{x},\mathbf{a};\boldsymbol{\theta})$, from which the components of the strain vector $\boldsymbol{\epsilon}$ are computed using automatic differentiation. Assume that the parameters are distributed according to a probability density given by $\rho(\mathbf{a})$ where $\mathbf{a} \in \omega$. A common problem in uncertainty quantification is to compute the first two statistical moments of a quantity of interest, which corresponds to computing the mean and variance (or standard deviation) of the lowest natural frequency in our case. These two statistics can be written as 

\begin{equation}\label{statistics}
    \mathbb{E}\Big( \lambda_1 \Big) = \int_{\omega} \rho(\mathbf{a}) \lambda_1 (\mathbf{a}) d\omega, \quad \text{Var}\Big(\lambda_1\Big)=\int_{\omega} \rho(\mathbf{a}) \lambda_1^2(\mathbf{a}) d\omega - \mathbb{E}\Big( \lambda_1 \Big)^2.
\end{equation}

In order to compute these quantities, we require the dependence of the eigenvalue on the parameters. Eq. \eqref{param1} could be solved for the eigenvalue as a function of the parameter one parameter setting at a time, where these values correspond to points on a grid discretizing the integrals in Eqs. \eqref{statistics}. However, the neural network discretization allows us to continuously approximate the eigenvalue's and eigenfunction's dependence on the parameters. The neural network is trained to construct the dependence of the eigenfunction on the spatial coordinate $\mathbf{x}$ and the material parameters $\mathbf{a}$ by minimizing the expected Rayleigh quotient

\begin{equation}\label{stochasticR}
    \underset{\boldsymbol{\theta}}{\text{argmin }} \int_{\omega} \qty(\frac{ \int_{\Omega} \boldsymbol{\epsilon} \cdot \mathbf{D}(\mathbf{x},\mathbf{a) \boldsymbol{\epsilon}} d\Omega }{\int_{\Omega} \mathbf{\hat u} \cdot \mathbf{\hat u}d\Omega})\rho(\mathbf{a}) d\omega.
\end{equation}

Just as we discretize integrals over the physical space $\Omega$, integrals over the space of stochastic parameters $\omega$ will also be discretized. Thus, Eq. \eqref{stochasticR} amounts to a sum of Rayleigh quotients at different parameter settings weighted by the probability density. Once the neural network discretizing the eigenfunctions is trained, the eigenvalue can be evaluated for specific parameter values at negligible cost.

 \paragraph{} To verify the implementation of the Rayleigh quotient in Eq. \eqref{param1}, we can compare it against an analytical solution for the case of a single material. A single material is obtained in the case where $a=0$ or $a=1$, which both correspond to the material interface being on one of the edges of the domain. We can choose the constitutive matrix $\mathbf{D}$ such that we solve the eigenvalue problem corresponding to

\begin{equation}\label{vectorlaplace}
    \nabla^2 \mathbf{u}(\mathbf{x}) + \lambda \mathbf{u}(\mathbf{x}) = \mathbf{0},
\end{equation}

\noindent for the given unit square domain and clamped boundaries. The exact solution for the first (normalized) eigenfunction is 

\begin{equation*}
    \mathbf{u}(\mathbf{x}) = \begin{bmatrix}
        2\sin \pi x_1 \sin \pi x_2 \\ 2\sin \pi x_1 \sin \pi x_2
    \end{bmatrix}.
\end{equation*}

The corresponding eigenvalue is $\lambda_1=2\pi^2=19.739$. Using a two hidden layer network with 10 hidden units in each layer and running 5000 epochs of ADAM optimization on Eq. \eqref{param1} with the appropriate constitutive matrix and a spatial integration grid with 2500 uniformly spaced points gives a converged Rayleigh quotient value of $\lambda_1=19.738$. This verifies the implementation of the deterministic Rayleigh quotient which shows up in the uncertainty quantification problem. 

\begin{figure}[hbt!]
\centering
\includegraphics[width=.95\textwidth]{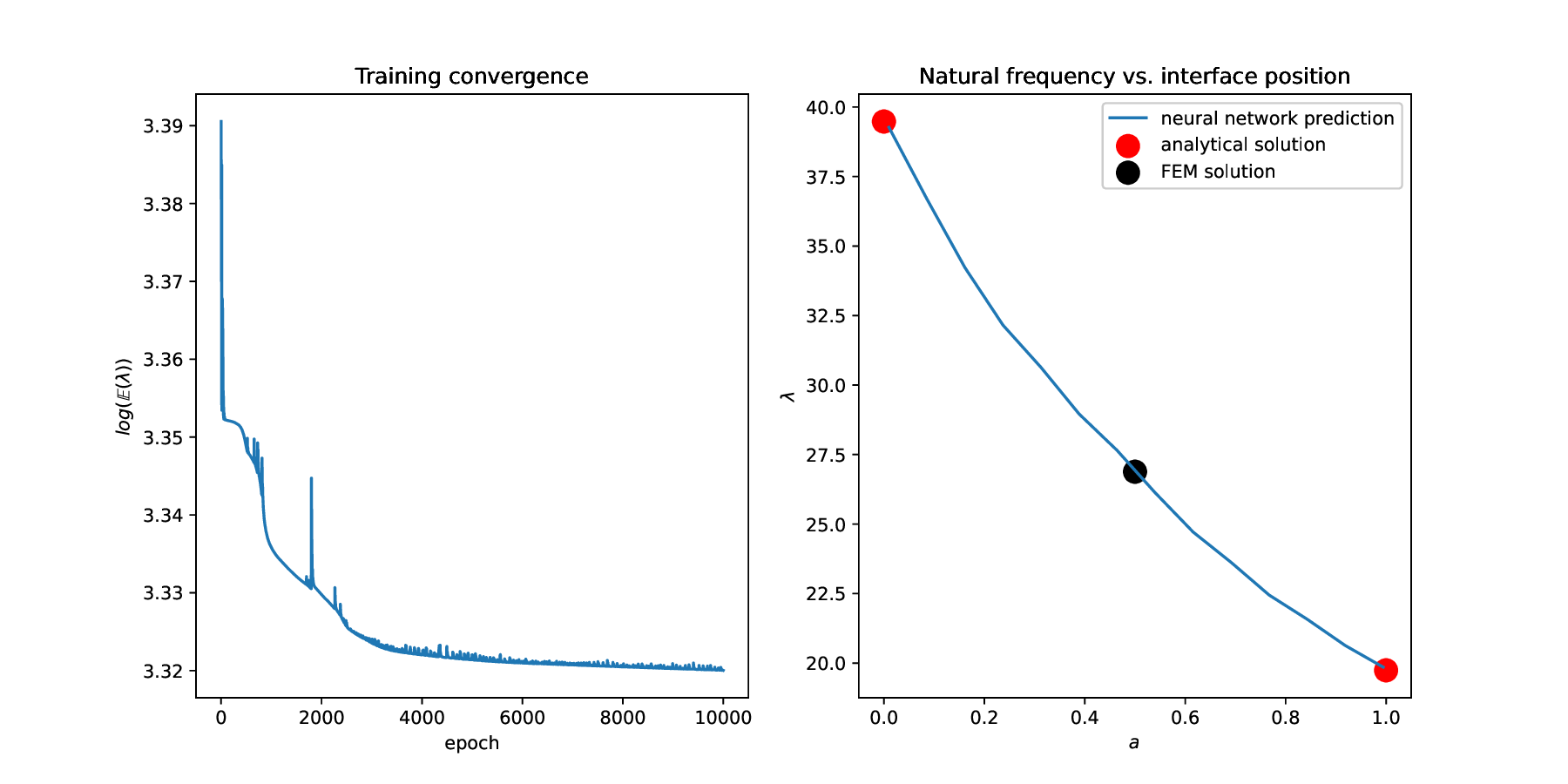}
\caption{Verifying the implementation of the expected Rayleigh quotient against analytical and FEM solutions for the first eigenvalue of the vector-valued Laplace equation. The stochastic parameter which controls the position of the interface between the two materials is included as an additional input to the neural network which discretizes the eigenfunction.}
\label{checklaplace}
\end{figure}

\paragraph{} As mentioned, we take the multimaterial modulus to be given by Eq. \eqref{material} with the interface position uniformly distributed as $a \sim \mathcal{U}(0,1)$. Taking Eq. \eqref{stochasticR} as the objective, we train the network to approximate the dependence of the eigenfunction on the spatial coordinate and material interface parameter. We also verify the implementation of the stochastic Rayleigh quotient with the vector-valued Laplace equation by checking that analytical solutions of the eigenvalue are returned in the extreme cases $a=0$ and $a=1$. These two cases correspond to constant coefficients in front of the Laplace operator in Eq. \eqref{vectorlaplace}. Taking $E_0=1$ and $E_1=2$, we minimize the expected Rayleigh quotient with a two hidden layer network in which both layers have a width of 15 and ADAM optimization with a learning rate of $1 \times 10^{-2}$. The Rayleigh quotient is integrated with a spatial integration grid of 5625 evenly spaced points. The trained neural network is then used to predict the dependence of the eigenvalue on the interface position $a$. We again find that removing the penalty which enforces normalization of the eigenfunction tends to lead to better results. There are 20 points in the integration grid for the stochastic space used to compute the expected Rayleigh quotient in Eq. \eqref{stochasticR}. This is adequate given that we expect low frequency, monotonic dependence of the natural frequency on the interface position. To verify the solution with FEM, we use triangular elements with a linear interpolation of the displacement components. The mesh consists of 1482 elements and we verify that the mesh is converged with the modulus given by Eq. \eqref{material} by refining to 1992 elements and noting that the computed natural frequency for $a=1/2$ changes by $<0.1\%$. The results are shown in Figure \ref{checklaplace}. We solve the eigenvalue problem by forming the mass and stiffness matrices computed with plane stress elements of unit thickness and the constitutive matrix corresponding to the vector-valued Laplace equation. The natural frequency is convex in the interface position $a$, and the analytical/FEM solutions are accurately returned.

\begin{figure}[hbt!]
\centering
\includegraphics[width=.95\textwidth]{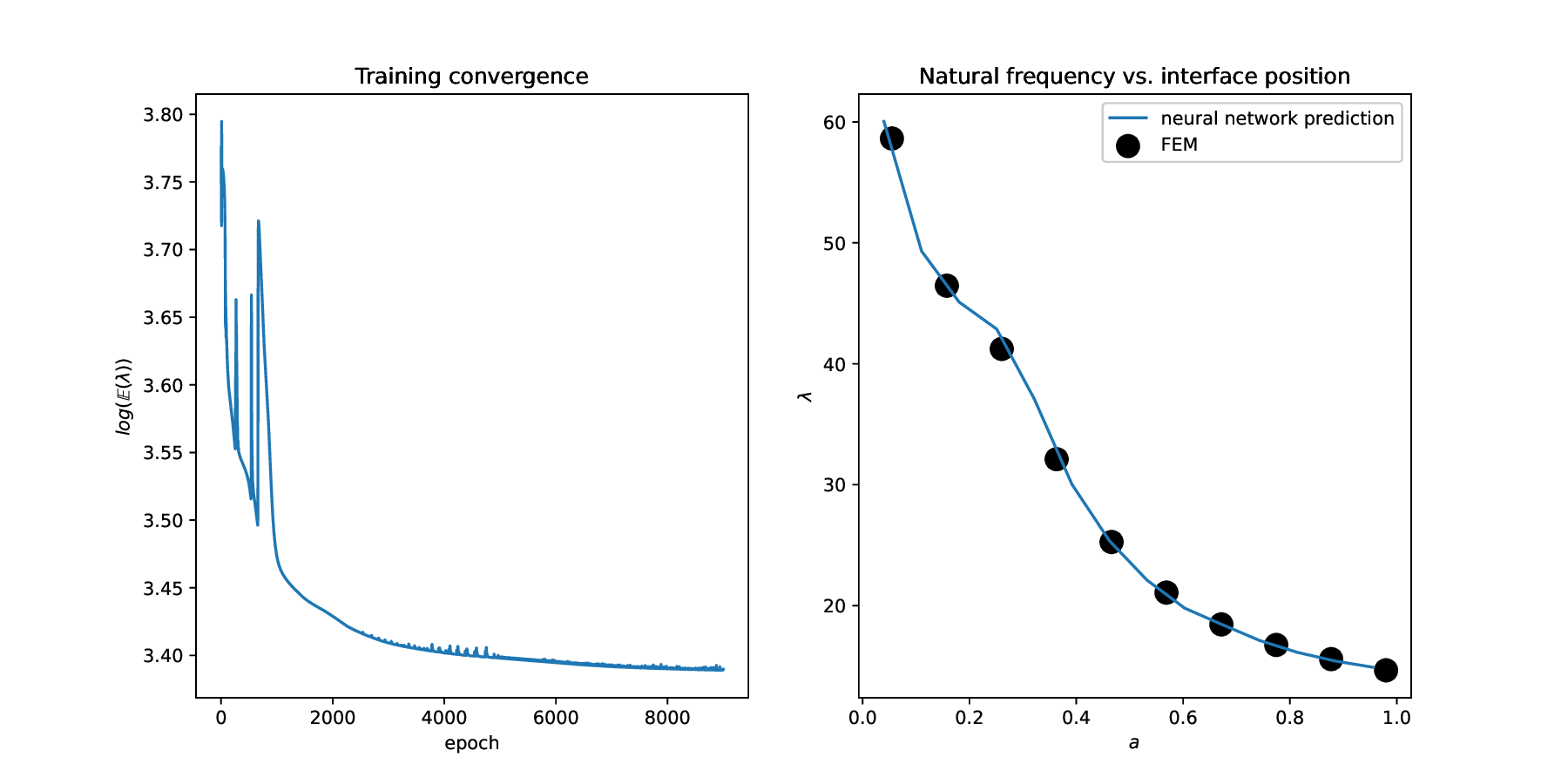}
\caption{The expected Rayleigh quotient is minimized to find the first natural frequency of the multimaterial structure as a function of the interface position. The predictions of the neural network agree with the finite element method results.}
\label{natfreq}
\end{figure}

\paragraph{} We now return to the problem of finding the natural frequency of the linearly elastic structure. We consider the configuration of Figure \ref{diagram} with clamped boundaries for both displacement components. The stochastic material interface parameter is uniformly distributed with a probability density given by $\rho(a) = \mathcal{U}(0,1)$. We will take the left and right modulus values to be $E_0=1$ and $E_1=5$ with a Poisson ratio of $\nu=0.25$. The plane stress constitutive matrix is given by Eq. \eqref{const}. We take 25 points in the integration grid for the stochastic space and use the same spatial integration grid of 5625 evenly spaced points. See Figure \ref{natfreq} for results. The predictions of the trained neural network for the natural frequency as a function of the interface position agree closely with the finite element results. Using the average relative error given in Eq. \eqref{toterror}, we have $\mathcal{E}_{\lambda}=1.4\times 10^{-2}$ when considering the 10 samples taken from FEM solutions. This example showcases the ability of neural networks to efficiently handle parametric problems, as there is no need to construct a separate basis to discretize the dependence of the solution on the parameter(s). Once the network is trained, it is computationally inexpensive to compute statistics of the quantity of interest. The mean and standard deviation are computed as 

\begin{equation*}
    \mathbb{E}(\lambda) = \int_0^1 \lambda(a) da = 29.42, \quad \sqrt{\text{Var}(\lambda)} = \sqrt{\int_0^1 \lambda(a)^2 da - \mathbb{E}(\lambda)^2 }= 14.25.
\end{equation*}

\subsection{Parametric geometry}

\begin{figure}[hbt!]
\centering
\includegraphics[width=0.7\textwidth]{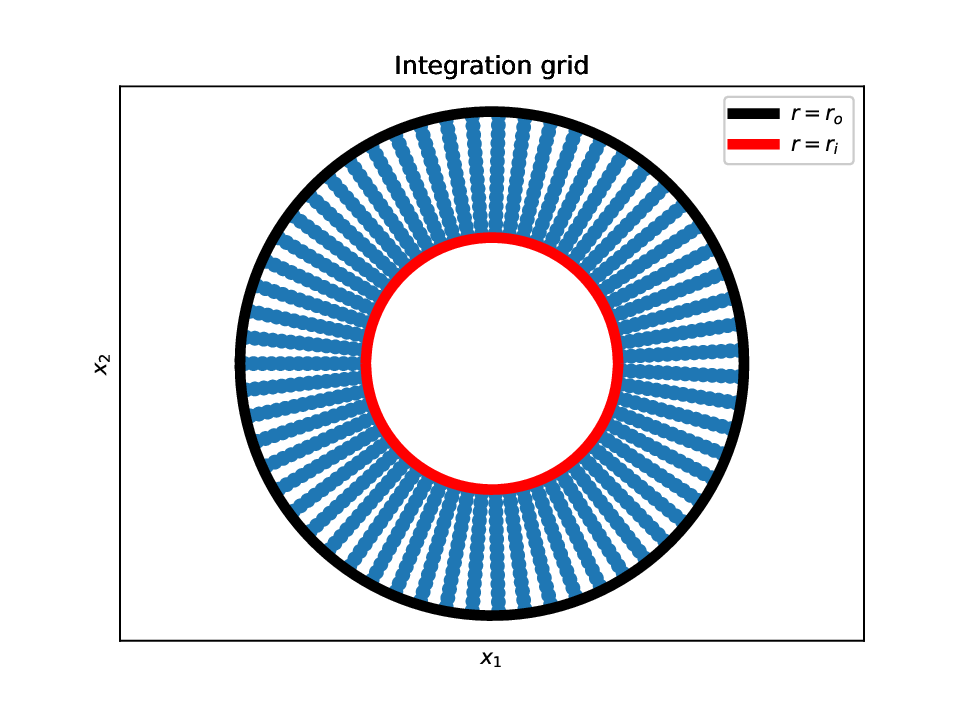}
\caption{Donut-shaped domain with given inner and outer radius. An example integration grid constructed in polar coordinates is shown in the figure. The integration weight corresponding to a point at radius $r$ is then $\Delta A = r \Delta r \Delta \theta$.}
\label{donutdiagram}
\end{figure}

\paragraph{} In addition to parameterizing the material, we can also parameterize the geometry of the computational domain $\Omega$. We will take the domain to be a donut-shaped region centered at the origin of the $x_1-x_2$ plane with fixed outer radius $r_o=1$ and a parameterized inner radius $r_i=a$. See Figure \ref{donutdiagram} for a schematic. We find a set of harmonic functions on the parameterized geometry by again working with the expected Rayleigh quotient. The neural network discretizing the dependence of the solution on the spatial coordinate $\mathbf{x}$ and the inner radius of the donut $a$ is trained with

\begin{equation*}
    \underset{\boldsymbol \theta}{\text{argmin }} \int_{\omega} \frac{\int \nabla \hat u (\mathbf{x},a;\boldsymbol \theta) \cdot \nabla \hat u(\mathbf{x},a;\boldsymbol{\theta}) d\Omega(a)}{\int u(\mathbf{x},a;\boldsymbol{\theta})^2 d\Omega(a)} \rho(a) d\omega.
\end{equation*}

We take the boundary conditions to be homogeneous and Dirichlet along both the inner and outer edges of the circle. This can be accomplished by multiplying the neural network output by

\begin{equation*}
D(\mathbf{x},a) = (1-x_1^2-x_2^2)(a^2-x_1^2-x_2^2).
\end{equation*}

The level-set function enforcing boundary conditions depends on the geometry parameter, as does the integration grid used to form the Rayleigh quotient. Because we look for a parametric set of eigenfunctions, we must modify the iterative algorithm given by Eq. \eqref{algorithm} to include dependence of the set of harmonic functions on the geometry parameter. The $i$-th optimization problem in the parametric setting is


\begin{equation}\label{algorithm2}
    \begin{aligned}
        \mathcal{L}^R_i = \mathbb{E}\Big[ \mathcal{R} \Big( \mathbf{\hat u}_i^{\perp}(\mathbf{x},a;\boldsymbol{\theta}) \Big) \Big] ,\\
        \mathbf{\hat u}_i^{\perp}(\mathbf{x},a;\boldsymbol{\theta}) = \mathbf{\tilde u}_i(\mathbf{x},a;\boldsymbol{\theta}) - \sum_{j=1}^{i-1} \frac{\int_{\Omega(a)}  \mathbf{\tilde u}_i (\mathbf{x},a;\boldsymbol{\theta})\cdot \mathbf{\hat u}_j(\mathbf{x},a) d\Omega(a)}{\int_{\Omega(a)} \mathbf{\hat u}_j (\mathbf{x},a) \cdot \mathbf{\hat u}_j(\mathbf{x},a) d\Omega(a)} \mathbf{\hat u}_j(\mathbf{x},a),\\
        \boldsymbol{\theta}_i = \underset{\boldsymbol{\theta}}{\text{argmin }} \mathcal{L}_i^R(\boldsymbol{\theta}).
    \end{aligned}
\end{equation}

\begin{figure}[hbt!]
\centering
\includegraphics[width=0.8\textwidth]{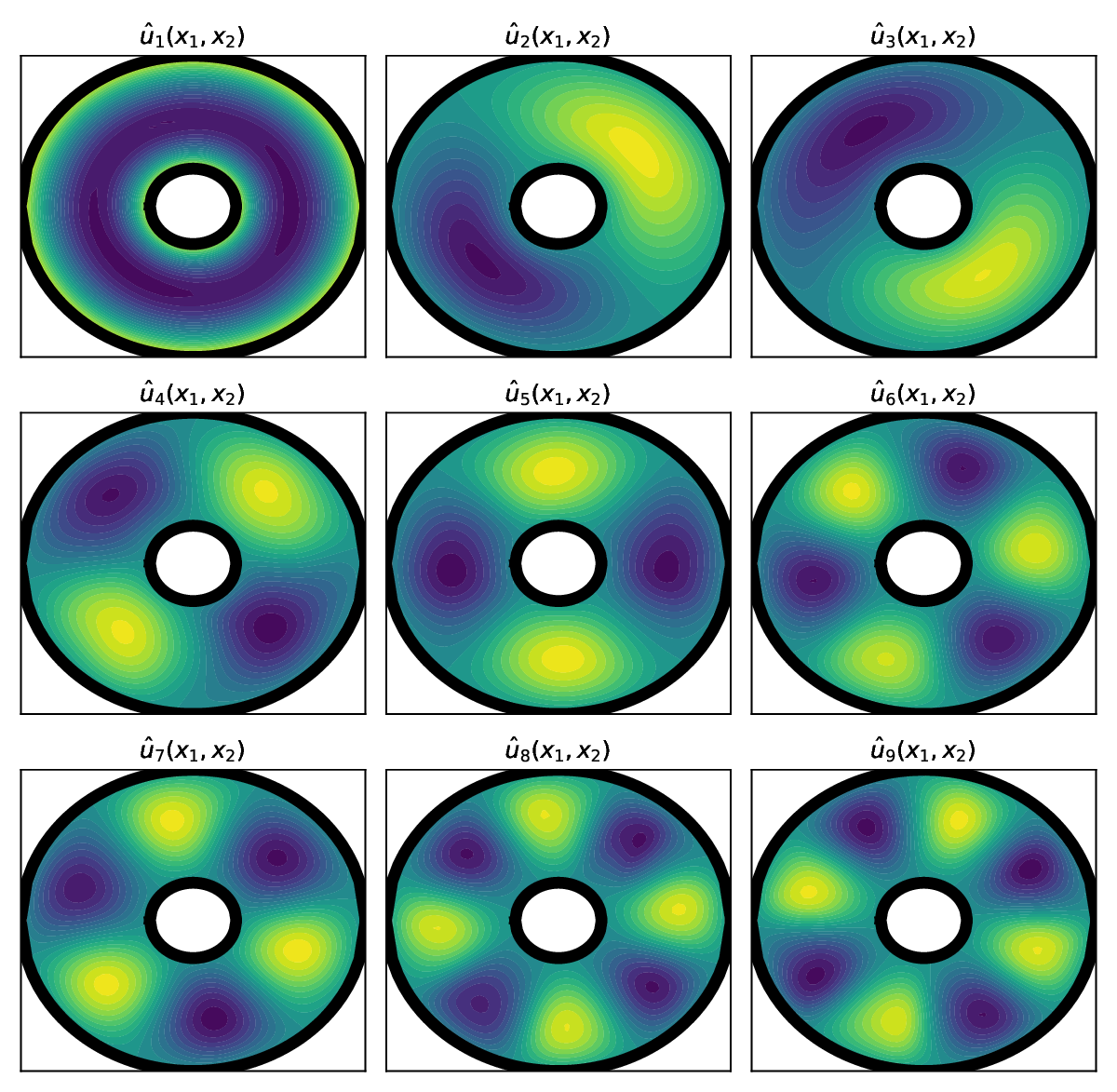}
\caption{The first 9 harmonic functions on the donut with homogeneous Dirichlet boundaries and inner radius $a=1/4$. We see that apart from the first, the eigenfunctions come in pairs that are 45 degree rotations of each other. The frequency increases with the index of the eigenfunction.}
\label{donut1}
\end{figure}

\begin{figure}[hbt!]
\centering
\includegraphics[width=0.8\textwidth]{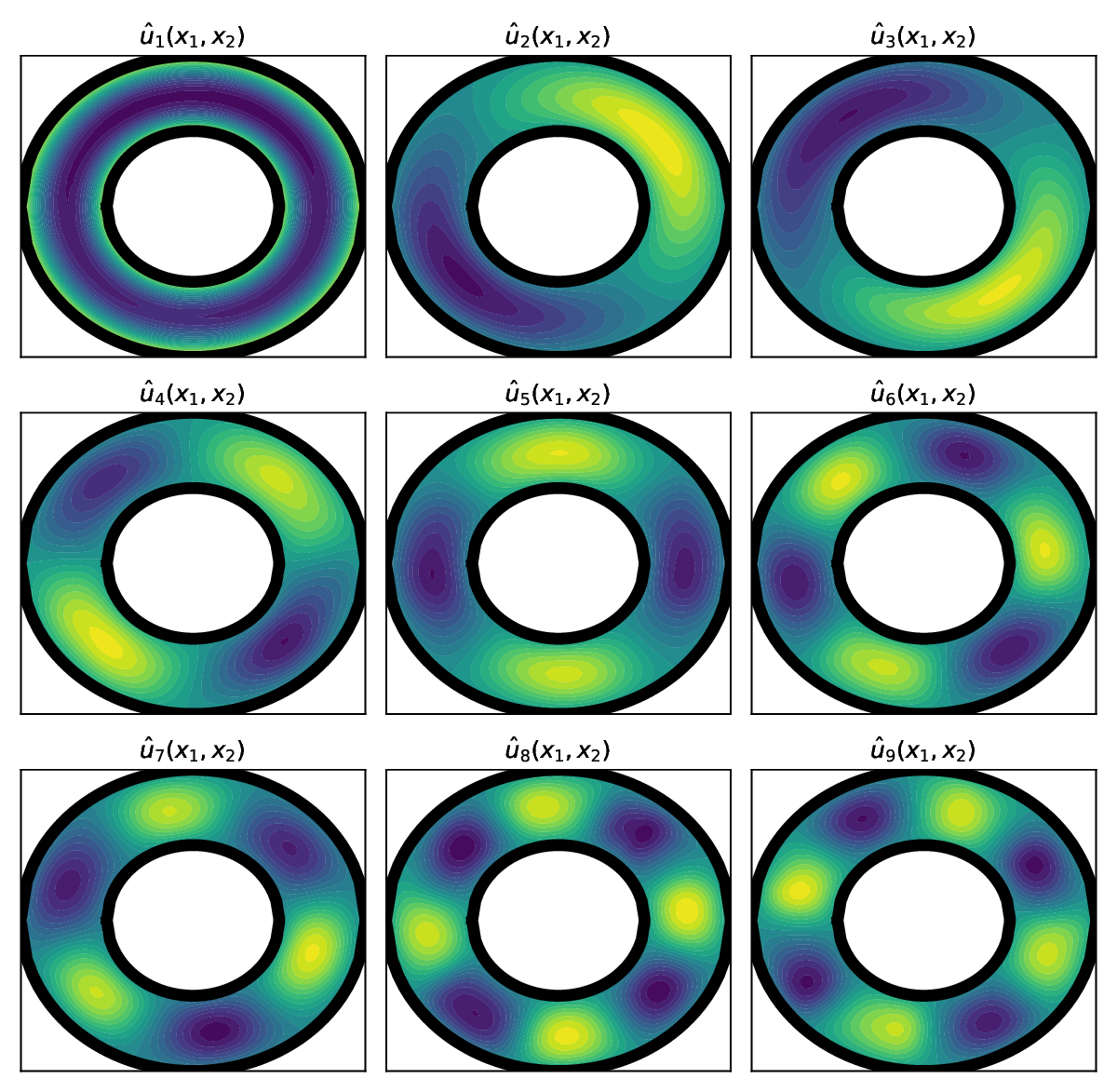}
\caption{The first 9 harmonic functions on the donut with homogeneous Dirichlet boundaries and inner radius $a=1/2$.}
\label{donut2}
\end{figure}

\begin{figure}[hbt!]
\centering
\includegraphics[width=0.75\textwidth]{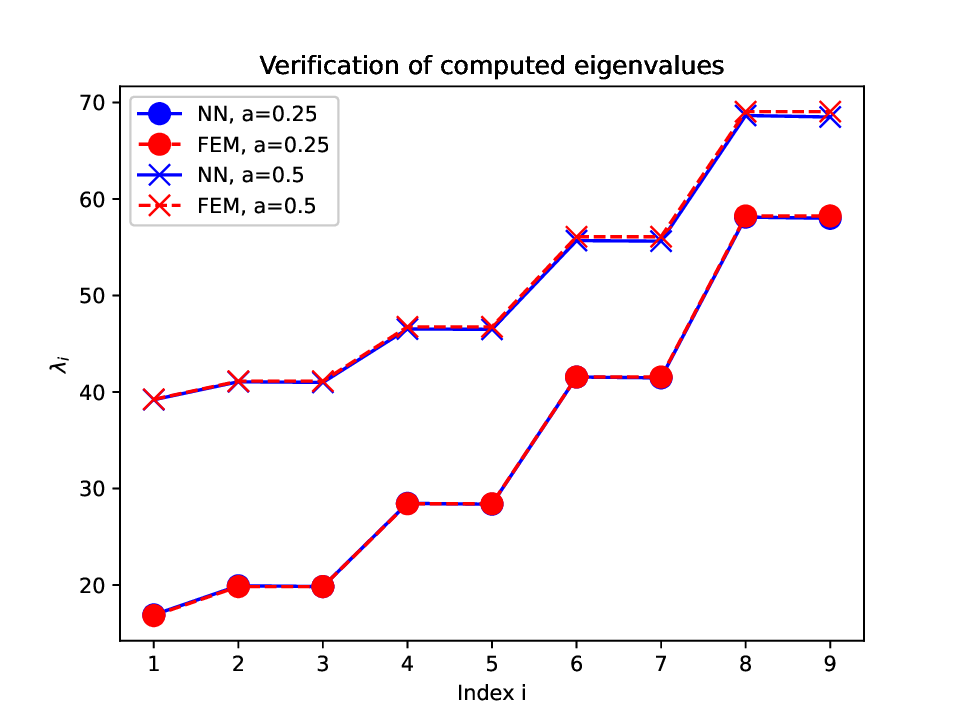}
\caption{The eigenvalues for each of the 9 harmonic functions with an inner radius of $a=1/4$ and $a=1/2$ are computed using FEM and the Rayleigh quotient method. We note that the Rayleigh quotient objective again provides ordered eigenpairs.}
\label{eigcomparison}
\end{figure}

In our implementation, we consider discrete points of the parameter space, and store integration grids for the domain $\Omega(a)$ corresponding to each one of them. We also store previously converged eigenfunctions at discrete settings of the parameter $a$ and use these in the orthogonalization. We note that this compact method to solve for sets of harmonic functions on parametric geometries has interesting applications in spectral methods. Basis functions corresponding to whole families of domain geometries can be constructed by solving the series of optimization problems given by Eq. \eqref{algorithm2}. Per the previous examples, ADAM optimization is run on each expected Rayleigh quotient until convergence is obtained. For simplicity, the parameter controlling the inner radius is taken to be a discrete random variable with $P(a=1/4) = P(a=1/2)=1/2$. Eq. \eqref{algorithm2} then computes the eigenpairs for two different geometries. The setup for the optimization problem is the same as before---we run the Rayleigh quotient optimization using ADAM optimization with a fixed learning of $1\times 10^{-2}$ for a set number of epochs once its value is within $\mathcal{T}=200$ of the previously converged Rayleigh quotient. The number of epochs after the convergence threshold is crossed is $5000+1000i$ where $i$ is the index of the eigenproblem. There are 6201 integration points for both geometric settings per the polar coordinate layout of Figure \ref{donutdiagram}. See Figures \ref{donut1} and \ref{donut2} for the first 9 harmonic functions at the two settings of the geometric parameter. For these settings, we again carry out eigenanalysis using the finite element using sufficiently refined meshes. See Figure \ref{eigcomparison} for a comparison of the eigenvalues computed by FEM and our Rayleigh quotient method. The average relative error (normalized by the FEM results) over the two parameter values is $\mathcal{E}_{\lambda} = 3\times 10^{-3}$. We note that it is more difficult to compare the eigenfunctions because they are only specified up to a rotational degree of freedom.

\section{Nonlinear eigenproblem}

\paragraph{} The p-Laplace equation is of mathematical interest as a nonlinear generalization of Laplace's equation. The scalar p-Laplace equation with a source term $s(\mathbf{x})$ is given by 
 
\begin{equation*}
    \nabla \cdot \Big( |\nabla u|^{p-2} \nabla u \Big) + s(\mathbf{x}) = 0 .
\end{equation*}

When $p=2$, the usual Laplace operator is recovered, leading to the scalar diffusion equation which governs steady state heat transfer. In the context of engineering mechanics, the p-Laplace equation can thus be thought of as a nonlinear generalization of heat transfer. For a given set of boundary conditions, the eigenvalue problem corresponding to the p-Laplace operator is

\begin{equation*}
    \nabla \cdot \Big( |\nabla u_i|^{p-2} \nabla u_i \Big) + \lambda_i |u_i|^{p-2} u_i = 0 .
\end{equation*}

Because of the nonlinearity, it is not possible to decouple the spatial and temporal components in a separation of variables solution to the transient heat conduction problem. This is the process typically used to arrive at the linear eigenvalue problems which provides the interpretation of eigenfunctions as spatial modes. Nonetheless, the Rayleigh quotient corresponding to the p-Laplace eigenproblem is

\begin{equation}\label{pr}
    \lambda = \underset{u(\mathbf{x})}{\text{min }} \frac{\int_{\Omega} |\nabla u|^{p} d\Omega}{\int_{\Omega} |u|^p d\Omega}.
\end{equation}

Physical interpretations aside, we will use the p-Laplace eigenproblem as a test case for nonlinear eigenanalysis using the proposed method. In two spatial dimensions, we can compare the performance of the neural network discretization of the Rayleigh quotient given by Eq. \eqref{pr} to a Fourier basis discretization in computing the first eigenvalue. The Fourier discretization for a unit square domain with homogeneous Dirichlet boundaries is given by 

\begin{equation*}\label{fourier}
\hat u^{\text{F}}(x_1,x_2;\boldsymbol{\theta}) = \sum_{i=1}^M\sum_{j=1}^M \theta_{ij} \sin( i \pi x_1 ) \sin( j \pi x_2),
\end{equation*}

\noindent where $M$ is the maximum frequency in each of the two coordinate directions. We take $M=10$ so there are $|\boldsymbol{\theta}|=100$ terms in the Fourier discretization of the solution. The neural network discretization is given by a two hidden layer neural network with a width of 7 in each of the two hidden layers, which corresponds to $|\boldsymbol{\theta}|=84$ trainable parameters. The neural network is initialized with Xavier initialization, and each coefficient in the Fourier series is initialized by drawing from a zero-mean normal distribution with standard deviation $\sigma=1 \times 10^{-3}$. The first eigenvalue is found by minimizing the Rayleigh quotient in Eq. \eqref{pr} with $p=5$ using full-batch ADAM optimization with a learning rate of $1\times 10^{-2}$ and an integration grid 5625 evenly spaced points. Like the Fourier series, the neural network satisfies the homogeneous Dirichlet boundaries by construction. Note that we again choose not enforce a penalty to normalize the eigenfunction for either the Fourier series or neural network discretization. 

\begin{figure}[hbt!]
\centering
\includegraphics[width=1.00\textwidth]{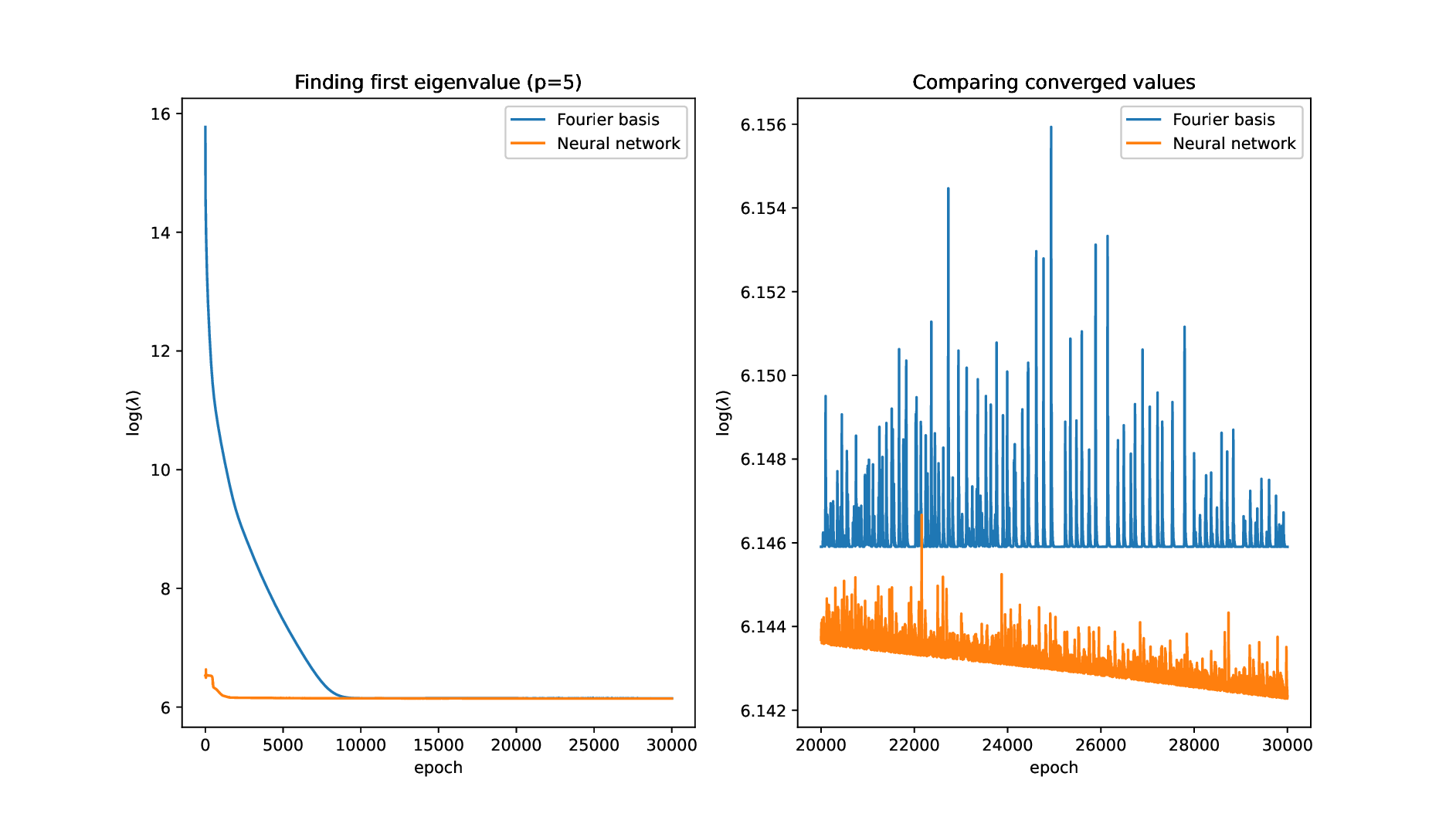}
\caption{Comparing the performance of a two hidden layer MLP neural network and Fourier basis for determining the first eigenvalue using the p-Laplace Rayleigh quotient. On the left, the whole training history is shown. On the right, we see a zoomed in picture of the Rayleigh quotient value near convergence. The neural network obtains a lower value, indicating higher accuracy.}
\label{plaplace}
\end{figure}

\paragraph{} See Figure \ref{plaplace} for a comparison of the solutions obtained from the two discretizations of the eigenfunction. The initial Rayleigh quotient for the Fourier series is extremely large, owing to improperly tuned high frequency contributions to the solution. This is an example where the spectral bias of the neural network, whereby low frequency solutions are preferred \cite{rahaman_spectral_2019}, is beneficial. The convergence of the neural network discretization is more rapid than the Fourier series. Although the neural network has fewer degrees of freedom than the Fourier series, the Rayleigh quotient after 30000 epochs is lower than the Fourier discretization, which has already saturated by this point. This indicates greater accuracy in both the eigenfunction and the eigenvalue. This is an interesting example where a neural network discretization outcompetes a more standard linear discretization in terms of computational cost. This is a consequence of the inherent nonlinearity of the Rayleigh quotient---which circumvents the need to form a linear system to solve the problem---and the spectral bias of the neural network representation of the solution.

\paragraph{} We note that it is simple to find a $p$-harmonic basis by modifying the Rayleigh quotient in Eq. \eqref{algorithm} to be the $p$-Laplace Rayleigh quotient. Numerical experiments indicate that learning the set of $p$-harmonic functions on a 2D square domain does not present any additional difficulties beyond those of the standard Laplace eigenproblem.

\section{High-dimensional eigenproblem}

\paragraph{} Another advantage of neural network approaches to solving eigenvalue problems with the Rayleigh quotient method is that it is easy to change the number of spatial dimensions. This contrasts with finite element discretizations of the solution field, which require significant redesigns of the basis functions if the dimension is changed. The Rayleigh quotient method also avoids matrix inversions in the solution of the eigenvalue problem, which becomes an asset when the dimensionality of the problem and the number of degrees of freedom in the solution discretization are large. We define the $d$-dimensional Laplace operator and the corresponding scalar eigenvalue problem as

\begin{equation}\label{dlaplace}
    \nabla^2_d := \sum_{i=1}^d \frac{\partial^2 }{\partial x_i^2}, \quad \nabla^2_d u(\mathbf{x}) + \lambda u(\mathbf{x}) = 0 .
\end{equation}

The first eigenvalue is the minimum of the $d$-dimensional Rayleigh quotient corresponding to Eq. \eqref{dlaplace}. This is given by 

\begin{equation*}
    \lambda_1 = \underset{u_1(\mathbf{x})}{\text{min }} \frac{\int_{\Omega} \nabla_d u_1 \cdot \nabla_d u_1 d\Omega }{\int_{\Omega} u_1^2 d\Omega},
\end{equation*}

\noindent where $\nabla_d$ is simply the $d$-dimensional gradient. We take the domain to be a unit hypercube $\Omega = [0,1]^d$ and the solution to have homogeneous Dirichlet boundaries enforced through Eq. \eqref{dirichlet}. Note that the exact solution for the first eigenfunction is given by 

\begin{equation*}\label{dd}
    u_1(x_1,\dots,x_d) = \prod_{i=1}^d \sin( \pi x_i) ,
\end{equation*}

with a corresponding eigenvalue of $\lambda_1=d\pi^2$. In the past examples, we have used sine functions to enforce Dirichlet boundaries on rectangular domains. Noting that the eigenfunction is a product of sines, we claim that this is an unfair (though coincidental) incorporation of prior knowledge into the solution for the eigenfunction. Thus, the Dirichlet boundaries will be enforced with a distance function of the form

\begin{equation*}
    D(x_1,\dots,x_d) = \prod_{i=1}^D x_i(1-x_i^4).
\end{equation*}

Due to its skew toward larger values of each coordinate $x_i$, this level set function avoids initializing the neural network discretization of the eigenfunction close to the exact solution. Though the neural network discretization is simple to modify to accommodate a $d$-dimensional spatial input, computing the integrals in the Rayleigh quotient poses challenges in high dimensions due to the curse of dimensionality. In \cite{yang_neural_2023}, high-dimensional Laplace eigenvalue problems are solved with wide and deep MLP neural networks by forming a large integration grid ($1 \times 10^5$ integration points) and training over a large number of epochs ($1 \times 10^5$ optimization steps). We argue that the convenience of neural network implementations should not be overshadowed by prohibitive computational expense. To this end, instead of generating and storing a fixed integration grid to compute the integrals, we will use Monte Carlo integration to form the Rayleigh quotient at each gradient step. This involves randomly sampling a set of integration points in the $d$-dimensional hypercube at each epoch of the training process, which is analogous to batching in stochastic gradient descent. Though Monte Carlo integration is often the preferred strategy in high dimensions, we note that, in order for this method to be effective, the decrease in the cost of each step must not lead to a corresponding increase in the number of steps until convergence. Writing the neural network discretization of the eigenfunction as $\hat u_1(\mathbf{x};\boldsymbol \theta)$ and treating the location of the integration points in each coordinate direction as independent random variables, the Rayleigh quotient is computed at each step with

\begin{equation}\label{sgd}
    \hat \lambda_1(\boldsymbol{\theta}) \approx \frac{\frac{1}{B}\sum_{i=1}^B \nabla_d \hat u_1(\mathbf{x}_i) \cdot \nabla_d \hat u_1(\mathbf{x}_i) }{\frac{1}{B}\sum_{i=1}^B \hat u_1(\mathbf{x}_i)^2},
\end{equation}

\noindent where $B$ is the number of integration points in the "batch" and each coordinate of $\mathbf{x}_i$ is drawn independently from a uniform distribution over the unit interval. Note that $B$ need not be so large that the integrals are accurately computed at each step, but this will introduce some noise into the gradients, which are also computed from Eq. \eqref{sgd} using automatic differentiation. By doing Monte Carlo integration, we avoid computing integrals on a $d$-dimensional grid at every step of the optimization, which becomes prohibitively expensive as $d$ becomes large. We run ADAM optimization to find the first eigenvalue for the hypercube by minimizing the Rayleigh quotient for dimensions $d=1,2,\dots,9$. We note that the architecture of neural networks makes it easy to implement a solution to this problem in an arbitrary number of dimensions. Knowing in advance that the first eigenfunction will have low frequency behavior, we use a two-hidden layer network, where each hidden layer has a width of 6. In our experience, wide and deep neural networks are unnecessary to accurately approximate low frequency functions. The Monte Carlo integration batch size is increased by 2000 for each successive dimension. The initial learning rate is $2 \times 10^{-3}$ and is multiplied by 0.95 every 1000 epochs.

\begin{figure}[hbt!]
\centering
\includegraphics[width=1.0\textwidth]{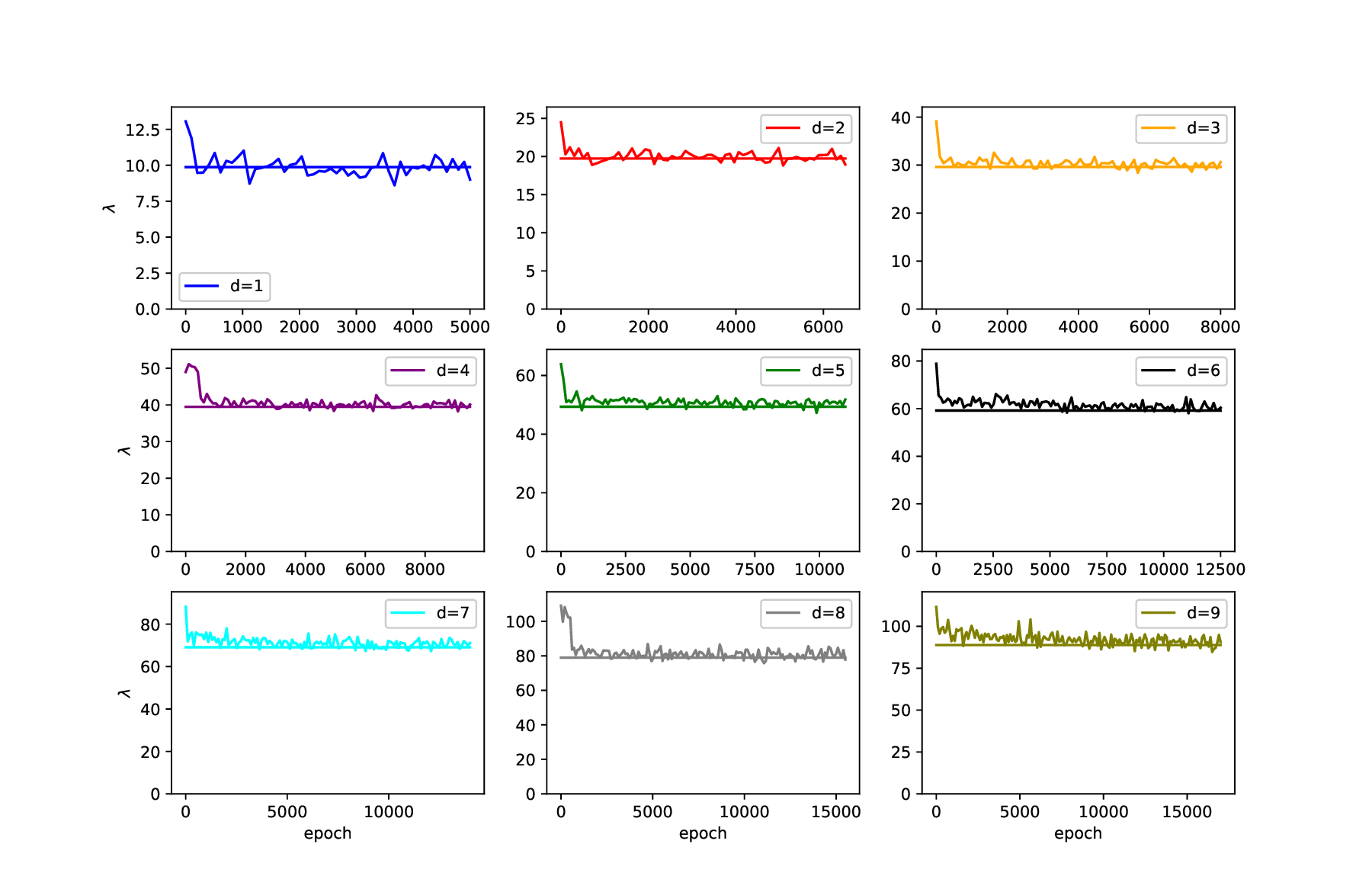}
\caption{Minimizing Rayleigh quotient with Monte Carlo integration for the Laplace eigenvalue problem in dimensions $d=1,2,\dots,9$. The true eigenvalue is shown as a solid line. The noise in the Rayleigh quotient arises from randomness in the value of the objective from Monte Carlo integration and the corresponding errors in the gradient computation.}
\label{highdim}
\end{figure}

\paragraph{} See Table \ref{tab:error} for the batch size and number of epochs at each dimension and Figure \ref{highdim} for the convergence of the Rayleigh quotient objective. We note that the Rayleigh quotient quickly approaches the true eigenvalue and then oscillates around it. If even a very sparse integration grid of 10 points per coordinate direction was used to deterministically solve this problem by forming stiffness and mass matrices with a more traditional spectral or finite element discretization, there would $10^9$ integration points in $d=9$ dimensions. Thus, even forming the 9-dimensional Rayleigh quotient for a one degree of freedom discretization with the unrealistically sparse integration grid of 10 points per coordinate direction requires more function evaluations than than our method, which uses a batch size of $1.7 \times 10^4$ integration points over $1.7 \times 10^4$ epochs. One might use a traditional finite element or spectral discretization with the Monte Carlo integration and iterative optimization strategy, but defining a traditional finite element mesh and basis in an arbitrary number of dimensions is a complex task. Furthermore, the number of degrees of freedom in both a finite element and spectral approximation scales poorly with the dimension. The number of nodes in a finite element mesh will scale like the integration grid and thus becomes intractably large. A spectral discretization has degrees of freedom that scale with $\binom{d+p}{p}$ where $p$ is the order of the discretization. Unless the eigenfunction happens to be an element of the spectral discretization, we require $p>1$ in order to obtain accurate results. Already at $p=3$, there are 220 degrees of freedom. In contrast, the neural network which discretizes the 9-dimensional Laplace eigenfunction has only 108 degrees of freedom and makes no assumptions about the functional form of the eigenfunction. This kind of parsimony should not be expected in general---high-dimensional functions with higher frequency behavior will require more trainable parameters to represent---but it is a striking feature of the neural network discretization that it is sometimes possible to represent high-dimensional functions with such few degrees of freedom. When the Rayleigh quotient is used to transform eigenanalysis into an optimization problem, neural networks are extremely convenient in the high-dimensional setting.

\paragraph{} We will estimate the eigenvalue by taking the Rayleigh quotient from the last 5000 epochs of the training and averaging it. The relative error with the true value is computed as 

\begin{equation}\label{errord}
    \mathcal{E}_{\lambda} = \frac{| \lambda - \langle \hat \lambda \rangle|}{\lambda},
\end{equation}

\noindent where $\langle \hat \lambda \rangle$ is the average Rayleigh quotient from the training and $\lambda$ is the true eigenvalue. See Table \ref{tab:error} for a summary of the results. As expected, the error increases with the dimension because the Rayleigh quotient and corresponding gradient are computed less accurately with the random integration scheme. However, averaging the Rayleigh quotient over the last 5000 epochs leads to small errors in the estimate of the eigenvalue even in high dimensions.

\begin{table}[h]
    \centering
    \begin{tabular}{|c|c|c|c|c|}
        \hline
        Dimension (d) & Integration batch size & Number of parameters & Epochs & Error (\%) \\ \hline
        1 & 1000 & 60 & 5000 & 0.6 \\ \hline
        2 & 3000 & 66 & 6500 & 1.0\\ \hline
        3 & 5000 & 72 & 8000 & 1.6 \\ \hline
       4 & 7000 & 78 & 9500 & 1.8 \\ \hline
        5 & 9000 & 84 & 11000& 1.7 \\ \hline
        6 & 11000 & 90 & 12500 & 2.0 \\ \hline
        7 & 13000 & 96 & 14000 & 2.1\\ \hline
        8 & 15000 & 102 & 15500 & 2.3 \\ \hline
        9 & 17000 & 108 & 17000 & 2.4 \\ \hline
    \end{tabular}
    \caption{Summary of results from the high-dimensional Laplace eigenproblem. }
    \label{tab:error}
\end{table}

\paragraph{} As a final study, we compute the second eigenvalue of the $d$-dimensional Laplace equation with the Monte Carlo integration scheme. As before, the first eigenvalue is computed by minimizing the Rayleigh quotient given by Eq. \eqref{sgd}. To compute the second eigenvalue, we need to orthogonalize the second candidate eigenfunction with respect to the converged value of the first eigenfunction. This requires integrals over the $d$-dimensional hypercube, which are intractable in the high-dimensional setting. We note that this is also the case when using sparse grid integration rules. Because we have an analytic form of the first eigenfunction, we can experiment with techniques for sparse grid integration in order to see the number of points required for accurate evaluation of the Rayleigh quotient \cite{gerstner_numerical_1998}. Even to compute $\int u_1^2 d\Omega$ with $u_1$ given by Eq. \ref{dd} in $d=10$ dimensions, we require on the order of $1 \times 10^6$ points using Smolyak sparse grid integration. This is a prohibitively large number of function evaluations in order to take a single optimization step. Thus, in addition to forming the Rayleigh quotient with Monte Carlo integration, we also perform Monte Carlo integration for the Gram-Schmidt orthogonalization procedure. The second eigenvalue is computed by minimizing

\begin{equation*}
    \hat \lambda_2 (\boldsymbol \theta) \approx \frac{\frac{1}{B} \sum_{i=1}^B \nabla_d \hat u^{\perp}_2(\mathbf{x}_i) \cdot \nabla_d \hat u^{\perp}_2(\mathbf{x}_i)}{\frac{1}{B} \sum_{i=1}^B \hat u_2^{\perp}(\mathbf{x}_i)^2}.
\end{equation*}

The Monte Carlo approximation of the Gram-Schmidt process for the second eigenfunction is given by 

\begin{equation*}
    \hat u_2^{\perp}(\mathbf{x}_i) \approx \tilde u_2(\mathbf{x}_i; \boldsymbol \theta) - \frac{\frac{1}{B}\sum_{i=1}^B \tilde u_2 (\mathbf{x}_i)\cdot \hat u_1(\mathbf{x}_i)}{\frac{1}{B}\sum_{i=1}^B \hat u_1(\mathbf{x_i})^2} \hat u_1(\mathbf{x_i}).
\end{equation*}

We use the same randomly sampled integration points for the Gram-Schmidt orthogonalization as are used in forming the Rayleigh quotient. Remember that $\hat u_1$ is the converged value of the first eigenfunction, and $\tilde u_2$ is the output of the neural network used to discretize the second eigenfunction. We will take $d=10$ and use the same hypercube with homogeneous Dirichlet boundaries, for which we have an exact solution for the second eigenvalue. As previously mentioned, the first eigenvalue is $\lambda_1=d\pi^2$. The second eigenvalue is given by $\lambda_1 = (d-1)\pi^2 + 4\pi^2$, which corresponds to an eigenfunction of the form

\begin{equation*}
    u_2(x_1,\dots,x_d) = \sin(2\pi x_1) \prod_{i=2}^d \sin(\pi x_i).
\end{equation*}

\begin{figure}[hbt!]
\centering
\includegraphics[width=1.0\textwidth]{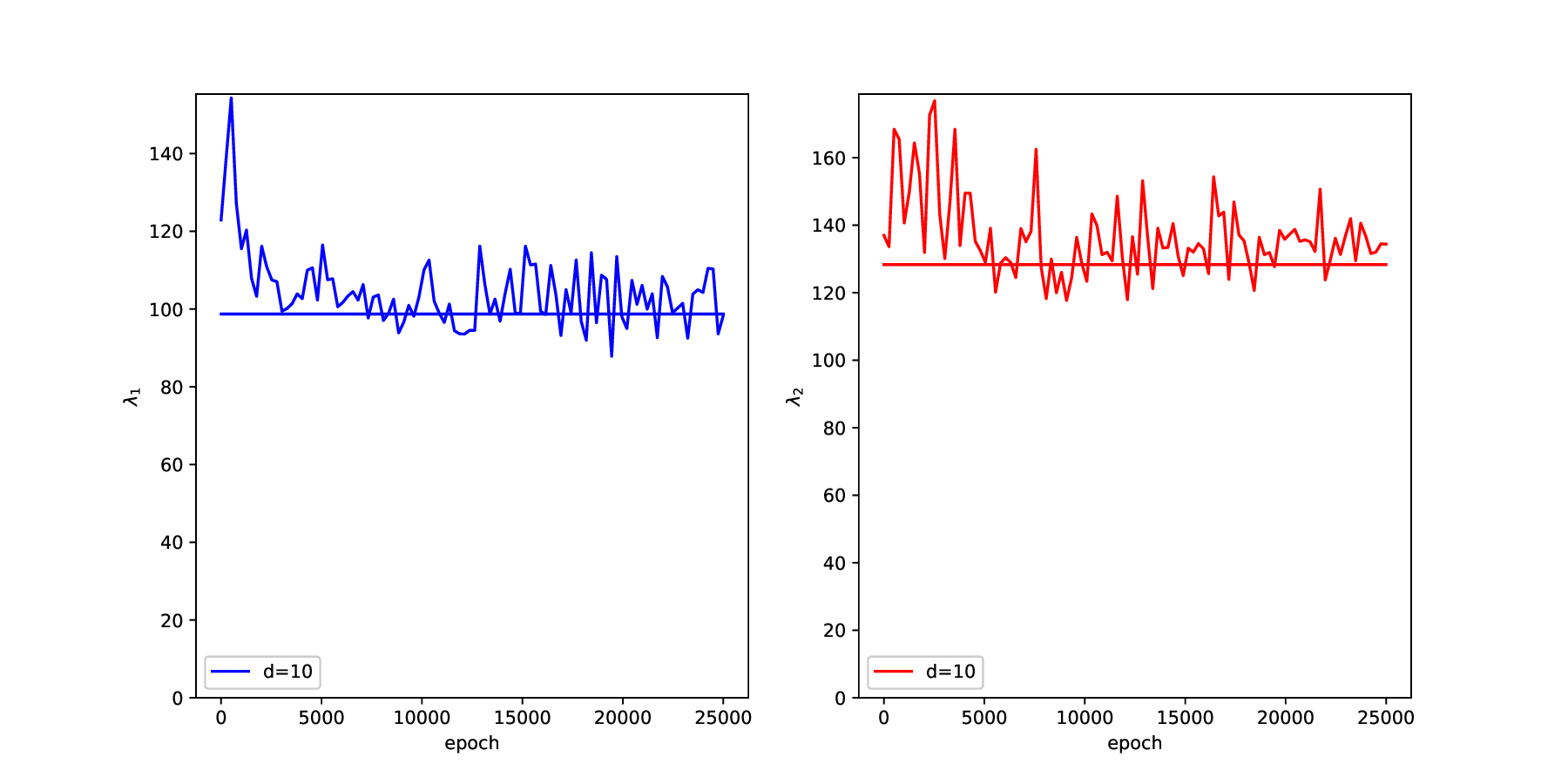}
\caption{Minimizing the Rayleigh quotient with Monte Carlo integration to estimate the first and second eigenvalue of the 10 dimensional Laplace equation. The eigenvalue is estimated by averaging the Rayleigh quotient value over the last 5000 iterations. The true eigenvalue is shown as a solid line.}
\label{d10}
\end{figure}

We take an integration batch size of $B=5000$ for the Rayleigh quotient minimization for both the first and second eigenvalue. The Rayleigh quotient minimization is cut off after 25000 epochs in both cases, and the eigenvalues are estimated by averaging the past 5000 values of the Rayleigh quotient. See Figure \ref{d10} for the results. Even though the Gram-Schmidt procedure introduces additional inaccuracies into the optimization problem for the second eigenvalue, we do recover a relatively accurate estimate. Using the error metric given in Eq. \eqref{errord}, the error of the first eigenvalue is $2.8\times 10^{-2}$, and the error of the second eigenvalue is $5.1\times 10^{-2}$. We note that in \cite{yang_neural_2023}, the first eigenvalue of the 10-dimensional Laplace equation is obtained with the power method. However, to compute the second eigenvalue with the power method, prior knowledge of its value is required to introduce an appropriate "shift" factor. Our method requires no knowledge of the spectrum of eigenvalues and remains computationally tractable. We see this as potentially compensating for the uncertainties in the final estimates of the eigenvalues introduced by the Monte Carlo integration.

\section{Conclusion}

\paragraph{} We have introduced a method to solve linear and nonlinear eigenvalue problems which uses neural network discretizations of the eigenfunctions, the Rayleigh quotient objective, and Gram-Schmidt orthogonalization to obtain an ordered sequence of eigenpairs. Though the combination of the Rayleigh quotient and orthogonalization has been used in the context of finite element methods before, it has not been explored in the physics-informed machine learning literature to the best of the authors' knowledge. The previous formulations of the eigenvalue problem in the PINNs literature rely on at least one of the following: 1) the strong form loss of the eigenvalue problem, which leads to non-unique solutions and does not return ordered eigenpairs, 2) penalties on the orthogonality of eigenfunctions, which leads to intractable optimization problems, or 3) the power method, which requires prior knowledge of the spectrum of eigenvalues. Our proposed method returns the eigenpairs in order, converges reliably, and requires no prior knowledge of the spectrum of eigenvalues. Another benefit of the Rayleigh quotient objective is that is naturally formulated as an optimization problem and it lowers the order of differentiation on the eigenfunction, which expedites computations involving automatic differentiation. Using our method, we have computed higher frequency eigenfunctions than any previous work in the PINNs literature. Furthermore, in showcasing our method, we highlighted the ability of neural networks to handle parametric and high-dimensional eigenvalue problems. For the reasons outlined above, we believe the combination of Rayleigh quotient objective and Gram-Schmidt orthogonalization to be well-suited for solving a wide range of eigenvalue problems discretized by neural networks.
 
\paragraph{} In future work, our method could be extended to solve parametric eigenproblems with a larger number of geometric and/or material parameters. Another potentially interesting area of future research is the intersection of machine learning techniques and the use of harmonic functions as spectral bases for numerical PDE solutions. Finding a set of harmonic functions on a family of parameterized geometries could be useful for quantifying uncertainty in the context of stochastic geometries. The ability of neural networks to interpolate is especially intriguing here, as harmonic functions could potentially be generated for an unseen geometric configuration, so long as the parameterized geometry is in the "interpolation" regime of the trained network. Regardless of the application, the crux of the proposed method of eigenanalysis is obtaining higher and higher frequency eigenfunctions. Even in one spatial dimension, eigenfunctions with a frequency of $n>14$ become very difficult to recover. This appears to be a combination of the inherent difficulties in fitting high frequency functions with neural networks (spectral bias) and numerical errors from the Gram-Schmidt process. In some sense, it is actually simpler to learn more eigenfunctions in higher dimensions, as the maximum frequency in each coordinate direction will always be less than the index of the eigenfunction. In other words, the frequency of the $i$-th eigenfunction in 1D is $i$, whereas it is less than $i$ in higher dimensions. From our experience, the difficulty of the optimization problem seems to increase with the maximum frequency of the target eigenfunction. Thus, future research should center around refining these methods such that very high frequency eigenfunctions can be reliably determined.


\bibliographystyle{plain}
\bibliography{My_Library}

\end{document}